%% file: main.tex
\documentclass{article}
\input{preamble}
\usepackage{quiver}
\usepackage[left=1.5in, right=1.5in]{geometry}
\usepackage{authblk}

\title{A Categorical Description of\\ Hierarchical Topic Modelling}
\author{Kaleb D. Ruscitti}
\affil{Department of Pure Mathematics, University of Waterloo}
\date{}

\newcommand{\Topic}{\mathrm{Topic}}

\newcommand{\Set}{\mathrm{Set}}
\newcommand{\Clust}{\mathrm{Clust}}

\begin{document}

\maketitle

\begin{abstract}
    This paper describes a hierarchical topic model of a dataset as a cluster functor from a poset of topics to the poset of fuzzy subsets of the dataset.  
    This presents a hierarchical topic model as a categorical database used to store the dataset.
    We construct a topic functor that maps fuzzy sets of data to their topics, and show it is a left adjoint to the clustering functor.
    We also describe data migration between the hierarchical topic model and databases storing metadata about the topics and samples.
    Finally we discuss implementation details and show examples of the database operations on real-world datasets. 
\end{abstract}

\tableofcontents

\paragraph{AI Contribution Statement}
Large language models (LLMs) were used as research and software development assistants in the preparation of this work. 
OpenAI's ChatGPT was used to review mathematical proofs for gaps or errors, discuss the computer implementation, and provide feedback on the organization of the manuscript.
Additionally, the GLM 5.3 and Qwen3.7 Plus models, accessed through OpenCode, were used to assist with writing and debugging the Python reference implementation.

The mathematical framework, research direction, initial manuscript draft, and scientific conclusions are the author's own work. LLM-generated suggestions and code were subject to review by the author, who takes full responsibility for the correctness of the mathematical arguments, software, and reported results.

\section{Introduction}

Topic modelling refers to unsupervised learning techniques that identify groups of related samples, called topics, and the relationships between them. Hierarchical topic modelling produces a \emph{topic tree}, where the children of a topic are its more specific subtopics \cite{griffiths_hierarchical_2003}. Topic models have proven useful in exploratory data analysis \cite{strydom_topic_2023, liang_molecules_2025, galli_topic_2025} and for content recommendation systems \cite{wilson_improving_2014, kawai_topic_2022}. 

Hierarchical topic models can be produced using embedding models to vectorize a dataset and then using hierarchical clustering techniques to obtain a topic tree \cite{grootendorst_bertopic_2022}. They can also be extracted from a flat classification of a dataset \cite{tibely_extracting_2013} or using formal concept analysis \cite{cimiano_learning_2005}. 
Many existing corpora have human-curated hierarchical topic models such as the Dewey Decimal System for library collections or the AMS Mathematical Subject Classifications for mathematical literature. 

Consider an organization with a large document corpus $X$ and associated metadata. 
An analyst can construct a hierarchical topic model by combining semantic embeddings, hierarchical clustering, and automated topic modelling tools. 
The resulting hierarchy can be equipped with an inclusion strength function that assigns to each document and cluster a membership value in $[0,1]$, with hard cluster assignments corresponding to the special case of binary membership values.
One particular hierarchical topic modelling procedure motivating this work is the Toponymy library \cite{healy_toponymy_2024} for constructing interpretable hierarchical topic models from unstructured text.
In working with these models, questions naturally arise about how topics should be represented mathematically, how collections of documents can be summarized by topics, and how topical information can be combined with other forms of data selection.

Although hierarchical clustering and topic modelling are useful and methods for constructing hierarchical topic models are well-studied, the mathematical structure of the operations performed on hierarchical topic models after they have been constructed has received less attention.
In particular, we are interested in a unified description of the relationship between topic descriptions and collections of samples, and of the interaction between topical operations and relational database queries.

The first goal of this paper is to provide a categorical database theory description of hierarchical topic models and their most important database operations.
This is the subject of Section \ref{s:hierarchical-clustering-defn}. 
This description views a hierarchical topic model as a functor $\Clust$ from the Lindenbaum-Tarski algebra associated to the tree $T$ into the poset of $[0,1]$-valued functions on the sample set $X$, under the domination ordering.
The Lindenbaum-Tarski algebra is a frame whose elements are (equivalence classes of) \emph{geometric formulas} built from the vertices of $T$ by applying conjunction and disjunction operations. 
The introduction of this algebra was initially motivated by the construction of an adjoint functor to $\Clust$, but the computational examples in Section \ref{sec:examples} demonstrate that this algebra intuitively describes the topics in a corpus and their relationships, and is a natural language to describe topic queries for a hierarchical topic model.
For the remainder of this paper, we consider a \emph{topic} to mean an equivalence class of geometric formulas. The vertices of the tree are the primitive topics from which all other topics are built.
These formulas are naturally human-interpretable and, combined with the functorial description of the operations of interest, provide a natural semantic interpretation of common operations applied to hierarchical topic models. 

In Theorem \ref{thm:topic-cluster-adjunction}, we prove that $\Clust$ has a left adjoint functor $\Topic$.
Semantically, $\Topic$ maps a $[0,1]$-valued function on $X$ (a fuzzy set of samples) to the most specific topic in the hierarchical clustering that contains the subset.    
Using the functors $\Topic$ and $\Clust$ we can prove results about operations used in applications of hierarchical topic models. For example, their composition $\Clust\circ\Topic$ is a closure operator (Corollary \ref{cor:closure}) that expands a subset of samples to the full set of samples on the same topic. This operation is useful for content recommendation, and one example of this application is described in Section \ref{sss:semantic-topical-integration}. 
 
The second goal is to understand how queries on hierarchical topic models interact with queries on sample metadata and topic metadata tables. 
Topic models are most useful when both the samples and the topics have rich metadata. 
This metadata can be stored in a relational database, and such databases already have a rich categorical description \cite{spivak_functorial_2012, spivak_relational_2015,schultz_algebraic_2017}.  
Oftentimes an analyst wants to compose queries between the relational and hierarchical database structures on the sample set. 
For example, the analyst might select documents written by a particular author and then identify the topics represented among them. Conversely, the analyst might select documents belonging to a particular topic and then investigate their metadata.
To compose hierarchical topic model query operations with relational database query operations, we need to understand the \emph{data migration functors} between these two different database structures on $X$. 
This is the subject of Section \ref{s:metadata-migrations}. 

The contribution of this paper is a unified compositional query framework for hierarchical topic models, together with sample and topic metadata. 
The framework extends the standard notion of topics as nodes in a cluster hierarchy to a distributive lattice of topic expressions that captures the composite topics that arise in real-world data. 
Data and topics are related through a pair of adjoint functors. 
The underlying lattice and adjunction constructions are standard; the contribution lies in using these constructions to organize topical queries and their interaction with metadata queries. 
The examples illustrate how composite topic expressions and metadata selections can be combined to formulate queries on document collections.

In addition to this theoretical model, there is a basic Python implementation, available on GitHub at the URL
\begin{center}
    \href{https://github.com/kalebruscitti/thematic-search}{https://github.com/kalebruscitti/thematic-search}.
\end{center}
An overview of the implementation details is provided in Section \ref{s:implementation}. Finally, Section \ref{sec:examples} describes examples of hierarchical database operations on two real-world datasets. Section \ref{ss:20newsgroups} looks at the 20 Newsgroups dataset, which comes with a human-curated topic model, and Section \ref{ss:ex-unsupervised} looks at a Wikipedia dataset with a machine-learned topic model generated with the Toponymy library.

\section{Categorical Description of Hierarchical Clustering}
\label{s:hierarchical-clustering-defn}

In this section, we will present a soft hierarchical clustering of a dataset $X$ as a functor $\Clust$ from a schema category $\mathcal{G}_T$ to fuzzy subsets of $X$.
Objects in the schema category are interpreted semantically as topics, and the category encodes the cluster tree hierarchy, built using a tree $T$. 
The cluster functor $\Clust$ records the inclusion strength of the samples in each cluster. 
After defining $\mathcal{G}_T$ and $\Clust$, we show that $\Clust$ has a left adjoint $\Topic$ (Theorem \ref{thm:topic-cluster-adjunction}), which takes fuzzy subsets of $X$ and returns the most specific topic in the hierarchical clustering that contains the subset. 

To make these ideas precise we will first need to review facts about posets, fuzzy sets, and geometric formulas. 

\subsection{Posets and Fuzzy Sets}

Let $P=(|P|,\leq)$ be a poset. A poset is a category with objects $|P|$ and with a single morphism $x\to y$ whenever $x\leq y$ in $P$. Given two posets $P$ and $Q$, a functor $P\to Q$ is equivalent to a monotone map, meaning a map $\phi:|P|\to |Q|$ such that if $x\leq y$ then $\phi(x)\leq \phi(y)$.

Let $X$ be a set. A \emph{fuzzy set contained in} $X$ is a map $f:X\to [0,1]$. Given two fuzzy sets $f$ and $g$, we say $f$ \emph{dominates} $g$, written $g\sqsubseteq f$ if $g(x)\leq f(x)$ for all $x\in X$. 
Let $[0,1]^X$ be the poset of fuzzy sets contained in $X$, under the domination ordering. 
We also define intersection and union of fuzzy sets: 
\begin{equation*}
    (f\cup g)(x) = \max(f(x),g(x)), \qquad (f\cap g)(x) = \min(f(x),g(x)).
\end{equation*}
These operations make $[0,1]^X$ into a distributive lattice \cite{zadeh_fuzzy_1965}. 
If a functor (of posets) between lattices also preserves the meet and join operations, then it is called a lattice homomorphism.
For any subset $U\subseteq X$, the characteristic function $\chi_U:X\to\{0,1\}$
\begin{equation*}
    \chi_U(x) = \begin{cases}
        1, & x\in U\\
        0, & x\not\in U
    \end{cases}
\end{equation*}
is a fuzzy set in $[0,1]^X$. Conversely if the range of a fuzzy set $f\in [0,1]^X$ is contained in $\{0,1\}$ then we call $f$ a set. 

\subsection{Geometric Formulas}
\label{ss:geometric-formula}

Suppose $P$ is a finite nonempty poset. 
A \emph{geometric formula} of $P$ is an expression of one of the following forms:
\begin{enumerate}
    \item (Atomic) $p$, for some $p\in P$.
    \item (Conjunction) $\psi \wedge \phi$, where $\psi$ and $\phi$ are geometric formulas.
    \item (Disjunction) $\psi \vee \phi$, where $\psi$ and $\phi$ are geometric formulas.
\end{enumerate}
A \emph{sequent} is an expression of the form $\psi_1,...,\psi_k \vdash \phi$, which should be read as ``if all of $\psi_1$ up to $\psi_k$ are true, then $\phi$ is true".
$P$ defines a set of sequents $S_P = \{p \vdash q ~:~ p\leq q\}$, and if a sequent $\phi \vdash \psi$ can be derived from $S_P$ and the rules of propositional logic, then we write $S_P \vdash \phi \vdash \psi$.
This defines a preorder on the set of geometric formulas, $\phi \leq \psi$ if $S_P \vdash \phi \vdash \psi$, and an equivalence relation by $\phi \sim \psi$ iff $\phi \leq \psi$ and $\psi \leq \phi$. 

Define $\mathcal{G}_P$ to be the set of equivalence classes of geometric formulas of $P$ under $\sim$, called the \emph{Lindenbaum-Tarski algebra} associated to $P$.
The preorder $\leq$ descends to a partial order on $\mathcal{G}_P$, and $\mathcal{G}_P$ is a finite distributive lattice with meet and join induced by $\wedge$ and $\vee$ \cite[Theorem 1.2]{dilworth_lattices_1945}.
As $P$ is finite, $\mathcal{G}_P$ has a minimal element $\bot :=\bigwedge_{P} p$ and a maximal element $\top := \bigvee_{P} p$.

To define lattice homomorphisms from $\mathcal{G}_P$ to other lattices $L$, in particular the lattice $[0,1]^X$, we will extend monotone maps from $P\to L$. 
The following lemma ensures that this is always well defined and unique.
\begin{lemma}
    \label{lem:unique-lattice-ext}
    Whenever $L$ is a distributive lattice and $f:P\to L$ is a monotone map, then there exists a unique extension of $f$ to a lattice homomorphism $\bar{f}:\mathcal{G}_P\to L$ such that for all $p,q\in P$,
\begin{enumerate}
    \item $\bar{f}([p]) = f(p)$,
    \item $\bar{f}([p\wedge q]) = f(p)\wedge f(q)$,
    \item $\bar{f}([p\vee q]) = f(p)\vee f(q)$.
\end{enumerate}
\end{lemma}
\begin{proof}
    This is essentially Proposition 4.4. of \cite{gehrke_vangool_2024}. 
    Define $\bar{f}$ recursively on geometric formulas by evaluating $\bar{f}([p]) = f(p)$ for all $p\in P$ and interpreting conjunctions and disjunctions as meets and joins. 
    We must verify that this evaluation depends only on the equivalence class of a formula.
    For a generating sequent $p\vdash q$
    \begin{equation*}
        \bar{f}([p]) = f(p) \leq f(q) = \bar{f}([q]).
    \end{equation*}
    Then as $\bar{f}$ is defined to preserve conjunction and disjunction and  $L$ is a distributive lattice, the logical inference rules for conjunction and disjunction preserve validity of $\bar{f}$ under this interpretation. 
    Then inducting on logical derivations lets one verify that
    \begin{equation*}
        S_P \vdash \phi \vdash \psi \implies \bar{f}(\phi) \leq \bar{f}(\psi).
    \end{equation*} 
    In particular if $\phi \sim \psi$ then $\bar{f}(\phi) =\bar{f}(\psi)$.
\end{proof}

When we use Lemma \ref{lem:unique-lattice-ext} we will abuse notation and denote the map from $P$ and from $\mathcal{G}_P$ with the same symbol.

\subsection{Hierarchical Clustering Functor}
\label{ss:soft-clustering}

Let $X\subset \mathbb{R}^d$ be a finite set. The data of a \emph{hierarchical clustering} of $X$ is:
\begin{itemize}
    \item A finite rooted tree $T$,
    \item For each vertex $v\in V(T)$, a fuzzy set $f_v:X\to [0,1]$.
\end{itemize}
If each $f_v$ is a set then we say the clustering is \emph{hard}, and otherwise it is \emph{soft}.
Let $V(T)$ be partially ordered by defining $v<w$ whenever $v$ is a descendant of $w$ and let $T$ also denote the poset $(V(T),\leq)$.
Under this partial order the root is the unique maximal element and we will denote it by $\top$.
Let us assume that the hierarchical clustering satisfies:
\begin{enumerate}
    \item $v \leq w \implies f_v \sqsubseteq f_w$.
    \item $f_{\top}(x) = 1$ for all $x \in X$.
\end{enumerate}
Under these assumptions, the map $\Clust:T\to [0,1]^X$ defined by $\Clust(v) = f_v$ for all $v\in V(T)$ is a functor. Let $\mathcal{G}_T$ be the Lindenbaum-Tarski algebra of $T$, defined in Subsection \ref{ss:geometric-formula}. We extend $\Clust$ to the unique lattice homomorphism $\Clust:\mathcal{G}_T\to [0,1]^X$ by defining
\begin{align*}
    \Clust(v\wedge w) &= \Clust(v)\cap \Clust(w), & \Clust(v\vee w) &= \Clust(v)\cup \Clust(w).
\end{align*}
As $T$ is finite, $\mathcal{G}_T$ is a finitely presented category and hence a valid schema category for categorical databases. 
In categorical database theory, a database is a functor from a schema category to $\Set$ \cite[Def. 3.2.8]{spivak_functorial_2012}. 
In this case we have a functor to fuzzy sets, so one might call $\Clust$ a `fuzzy categorical database'. 

\begin{figure}
    \centering
    \includegraphics[width=0.5\linewidth]{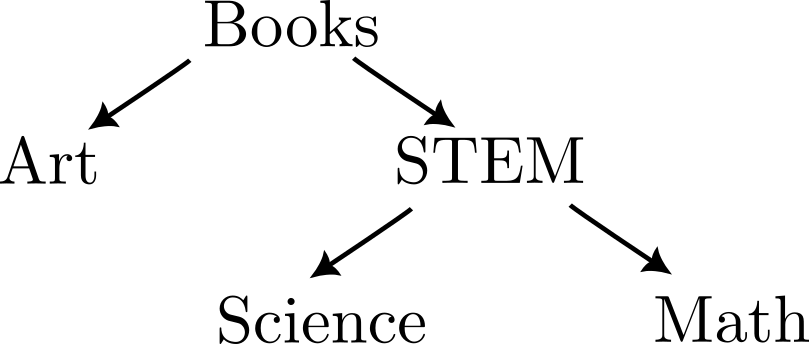}
    \caption{A small example topic tree for a dataset of books. The root of this tree is 'Books'.}
    \label{fig:small-topic-tree-example}
\end{figure}

The extension from $T$ to $\mathcal{G}_T$ is motivated both theoretically and from applications. 
Theoretically the extension is required to make Theorem \ref{thm:topic-cluster-adjunction} true, but the non-atomic geometric formulas in $\mathcal{G}_T$ also arise naturally in applications. 
To see why, let's consider the example topic tree in Figure \ref{fig:small-topic-tree-example}. 
How might the geometric formulas `Art $\wedge$ Math' or `Art $\vee$ Math' arise, and how do they differ? 
Suppose this tree was generated for a collection $X$ of books and $u_i:X\to \{0,1\}$ is a subset of $X$ for $i=1,2$, perhaps representing the expressed preferences of two people. 
Consider 
\begin{align*}
    u_1 &= \{`\textit{Painting for Dummies}', `\textit{Categories for the Working Mathematician}'\},\\ 
    u_2 &= \{`\textit{G\"odel, Escher, Bach}'\}.
\end{align*}
The subset $u_1$ contains one book about art, and another about math. 
Thus it is reasonable to say that the topic of $u_1$ is 'Art $\vee$ Math'.
On the other hand, \textit{G\"odel, Escher, Bach} (by Douglas Hofstadter) is a book discussing similarities in the work of the logician G\"odel, illustrator Escher, and composer Bach, and so the topic of $u_2$ would be better described as `Art $\wedge$ Math'. 
The goal of the next subsection is to make the topic of a fuzzy subset of $X$ precise.

\subsection{The Topic Functor}
\label{ss:topic-functor}

Given a hierarchical clustering $\Clust:\mathcal{G}_T\to [0,1]^X$ of a dataset $X$ and a fuzzy subset $f:X\to [0,1]$, we now define a functor sending $f$ to the most specific topic whose cluster contains $f$.
\pagebreak 

\begin{theorem}
\label{thm:topic-cluster-adjunction}
Suppose $\Clust:\mathcal{G}_T\to [0,1]^X$ is a hierarchical clustering. Then the functor $\Clust: \mathcal{G}_T \to [0,1]^X$ has a left adjoint $\Topic: [0,1]^X \to \mathcal{G}_T$ defined by:
\begin{equation}
    \Topic(f) = \bigwedge\{\psi\in \mathcal{G}_T ~|~ f \sqsubseteq \Clust(\psi)\}
\end{equation}
\end{theorem}
\begin{proof}
The set being met is non-empty because $f\sqsubseteq \Clust(\top)$. 
The meet always exists because $\mathcal{G}_T$ is closed under finite meets by construction, and as $T$ is finite, so $\mathcal{G}_T$ is also finite.
Being a left adjoint means that for all $f\in [0,1]^X$ and $\phi\in \mathcal{G}_T$:
$$\Topic(f) \leq \phi \iff f \sqsubseteq 
\Clust(\phi).$$
($\Leftarrow$) If $f \sqsubseteq \Clust(\phi)$, then $\phi$ is a member of the set $\{ \psi \mid f \sqsubseteq f_\psi \}$. Therefore by the definitions of $\Topic(f)$ and $\leq$, $\Topic(f) \leq \phi$.\\
($\Rightarrow$) Suppose $\Topic(f) \leq \phi$. By the monotonicity of $\Clust$, 
$$\Clust(\Topic(f)) \sqsubseteq \Clust(\phi).$$ 
Moreover as $\Clust$ preserves meets $\Clust(\Topic(f))$ is a meet over fuzzy sets dominating $f$, and hence $f\sqsubseteq \Clust(\Topic(f))$. 
Combining these inequalities we obtain $f\sqsubseteq \Clust(\phi)$.
\end{proof}

A pair of adjoint functors between posets is called a \emph{Galois connection} in the field of formal concept analysis \cite{ganter_formal_1999,fong_invitation_2019}. 

\begin{corollary}
    The functor $\Topic:[0,1]^X\to \mathcal{G}_T$ preserves colimits, and the functor $\Clust:\mathcal{G}_T\to [0,1]^X$ preserves limits.
\end{corollary}
\begin{proof}
    As $\Topic\dashv \Clust$, this follows from the standard result that left and right adjoints preserve colimits and limits respectively \cite[\S V.5]{mac_lane_categories_1994}.
\end{proof}

One specific implication of this corollary is that for any two fuzzy sets $f,g\in [0,1]^X$, 
\begin{equation*}
    \Topic(f\cup g) = \Topic(f)\vee \Topic(g).
\end{equation*}
This is because unions are colimits in the category of fuzzy sets.
Another corollary of the adjunction $\Topic\dashv \Clust$ is that $\Clust\circ \Topic$ is a closure operator:
\begin{corollary}
    \label{cor:closure}
    Let $\mathcal{C}:=\Clust\circ \Topic:[0,1]^X\to[0,1]^X$. Then for all $f\in [0,1]^X$,
    $$\mathcal{C}(\mathcal{C}(f)) = \mathcal{C}(f)\quad \text{ and }\quad f\sqsubseteq \mathcal{C}(f).$$
\end{corollary}
\begin{proof}
    See \cite[\S0.4 Proposition 8]{ganter_formal_1999}
\end{proof}
Given any fuzzy set of samples $f\in [0,1]^X$, we call $\mathcal{C}(f)$ its \emph{topical closure}. This is the smallest fuzzy set of documents that both have a shared topic and also contain $f$.

\subsection{Topical Coarse-graining}
\label{ss:coarse-graining}

The tree structure of a hierarchical topic model $\Clust:T\to [0,1]^X$ imposes a natural notion of `coarse-graining' on fuzzy sets $f:X\to [0,1]$, in particular on the clusters of the topic model.

For $v\neq \top$ in $V(T)$, let $\text{Parent}(v) := \min\{w\in V(T) ~|~ v< w\}$, and let $\text{Parent}(\top):=\top$. 
This is a functor $T\to T$ because whenever $v< w$ then either $w=\text{Parent}(v)$ or $w$ is an ancestor of $\text{Parent}(v)$; in either case $\text{Parent}(v) \leq w \leq \text{Parent}(w)$.
Therefore there is a unique extension to a lattice homomorphism $\text{Parent}:\mathcal{G}_T\to\mathcal{G}_T$ by defining that Parent commutes with the binary operators $\wedge$ and $\vee$:
\begin{align*}
    \text{Parent}(v\wedge w) &= \text{Parent}(v)\wedge \text{Parent}(w), &     \text{Parent}(v\vee w) &= \text{Parent}(v)\vee \text{Parent}(w).
\end{align*}
By definition, $v \leq \text{Parent}(v)$ for every $v\in V(T)$. 
For a general element $\theta \in \mathcal{G}_T$ the inequality $\theta \leq \text{Parent}(\theta)$ can be proven by induction on the number of $\wedge$ and $\vee$ needed to express $\theta$. 
The inductive step is that if $\theta = [\phi \wedge \psi]$ then 
\begin{equation*}
    \phi \wedge \psi \leq \text{Parent}(\phi)\wedge\text{Parent}(\psi) = \text{Parent}(\phi \wedge \psi),
\end{equation*}
and the same holds replacing $\wedge$ with $\vee$.

\begin{definition}
    Define the \emph{coarse-graining} functor $\mathrm{Coarse}:[0,1]^X\to [0,1]^X$ to be
    $$\mathrm{Coarse}:= \Clust\circ \text{Parent}\circ\Topic.$$
\end{definition}

As the name suggests, coarse-graining a fuzzy set produces a larger fuzzy set:

\begin{proposition}
    \label{prop:coarse-graining}
    For any $f\in [0,1]^X$, $f\sqsubseteq\mathrm{Coarse}(f)$.
\end{proposition}
\begin{proof}
    Using that $\phi\leq \text{Parent}(\phi)$ for any $\phi\in \mathcal{G}_T$ and functoriality of the hierarchical clustering, we know $\Clust(\phi) \sqsubseteq \Clust(\mathrm{Parent}(\phi))$. Applying this to $\phi=\Topic(f)$ and using Corollary \ref{cor:closure} we get 
    $$f \sqsubseteq \Clust(\Topic(f)) \sqsubseteq \mathrm{Coarse}(f).$$
\end{proof}

  \subsection{Robust Topic Estimator}
    \label{ss:robust-topic}
    
    In Subsection \ref{ss:soft-clustering} we defined the functor $\Topic:[0,1]^X\to \mathcal{G}_T$ that maps a fuzzy set $f$ of samples in $X$ to the most specific topic that dominates the set.
    The strict requirement that $\Topic(f)$ dominates $f$ means $\Topic$ is very sensitive to small changes in $f$. 
    Suppose, for example, that $f$ contains 100 samples from topic $A$ with inclusion strength $1$, and a single sample from topic $B$ with inclusion strength $\epsilon>0$ that is not in topic $A$. 
    Then $A$ does not dominate $f$, so $\Topic(f)\neq A$. It could be instead that $\Topic(f) = A\vee B$.
    
    For real-world datasets with topic models generated by noisy statistical processes (e.g. unsupervised machine learning) it is likely that the values of $\Clust$ will be noisy and the sensitivity of $\Topic$ to this noise degrades the quality of $\Topic(f)$ as a topical description of the set $f$. 
    For this reason, we want to define a more robust version of $\Topic$.
    This is analogous  to choosing a parametric statistical model to fit a dataset:
    It is well known that there is a trade-off between model fit and complexity, and there are heuristics like the Bayesian information criterion (BIC) that are used to choose an appropriate statistical model \cite{schwarz_estimating_1978}.

    There are many ways that one could choose a robust topic for a set $f$. 
    Here we describe one reasonable approach inspired by the BIC. The BIC says that one should select a statistical model that maximizes $2\ell - k\log(n)$, where $\ell$ is the log-likelihood of the model, $k$ is the number of model parameters and $n$ is the number of samples. 
    Inspired by the BIC, we want to choose a topic that maximizes a similar score function, where we have replaced $\ell$ with a suitable metric of topic fit and $k$ with a suitable measure of the complexity of a topic.

    Let $f\neq0$ be a fuzzy subset of $X$ and let $|f| = \sum_{x\in X} f(x)$.
    In the definition of $\Topic$ we take the conjunction of all topics $\phi$ for which $f \sqsubseteq \Clust(\phi)$, meaning we are requiring that $\Clust(\phi)$ \emph{covers} the set $f$. 
    Motivated by this, we define the \emph{coverage score}
    \begin{equation*}
        c(f, \phi) = \frac{1}{|f|}\sum_{x\in X}\min\{\Clust(\phi)(x), f(x)\}.
    \end{equation*} 
    When $f \sqsubseteq \Clust(\phi)$, $c(f,\phi) = 1$. Otherwise $0\leq c(f,\phi) < 1$.
    Coverage is not sufficient to score topics however, as we want the \emph{most specific} topic that has good coverage of the fuzzy set. 
    For any topic $\phi\in\mathcal{G}_T$, let 
    $$p(\phi) = \frac{1}{|X|}\sum_{x\in X}\Clust(\phi)(x),$$
    which is the normalized size of the topic $\phi$. Then let $S$ be the \emph{surprise score}, modelled after Kullback-Leibler (KL) divergence:
    \begin{equation*}
        S(f,\phi) = c(f, \phi) \log \frac{c(f,\phi)}{p(\phi)}.
    \end{equation*}
    If $c(f,\phi)=0$ then we define $S(f,\phi)=0$. When $p(\phi)=0$ it must also occur that $c(f,\phi)=0$, so again $S(f,\phi)=0$.
    As $\Clust$ is monotone so is $p$, and hence whenever $\phi \leq \psi$, we know that $p(\phi) \leq p(\psi)$. This means that when $\phi\leq \psi$ have equal coverage of $f$, $\phi$ will be a smaller topic, and therefore have higher or equal surprise score. 
    Maximizing the score function is therefore a soft version of $\Topic$, so we take this surprise score as a measure of how well a topic $\phi$ fits a fuzzy set $f$ of samples from $X$. 

    Next we introduce a complexity penalty on elements of $\mathcal{G}_T$. 
    The complexity of a geometric formula $\phi$ can be described using a rank function $r$ defined recursively by $r(t) = 0$ for $t\in V(T)$ and 
    $$r(\phi \vee \psi) = r(\phi\wedge \psi) := r(\phi) + r(\psi) + 1.$$
    However this does not descend to a well defined function on $\mathcal{G}_T$. If it did, then because $\phi\wedge \phi = \phi$ we must have $r(\phi) = 2r(\phi)+1$, a contradiction.
    Instead, we can use the 'minimum description rank' of a class $\theta \in \mathcal{G}_T$:
    \begin{equation*}
        \hat{r}(\theta) := \min_{\phi: [\phi]=\theta} r(\phi). 
    \end{equation*}
    For any $\theta \in \mathcal{G}_T$, $\hat{r}(\theta)$ is the minimum number of $\wedge$ and $\vee$ operators needed to write $\theta$ in terms of the vertices of $T$.

    Then we finally define the robust topic function $\text{RobustTopic}:[0,1]^X\to \mathcal{G}_T$ as
    \begin{equation}
        \label{eqn:topic-score-function}
        \text{RobustTopic}(f) = \text{argmax}_{\phi\in\mathcal{G}_T} \left\{S(f,\phi) - \lambda \frac{\log(1+|f|)}{|f|} \hat{r}(\phi) \right\}, \qquad \lambda>0.
    \end{equation}
    If $f=0$ we declare $\text{RobustTopic}(f)=\top$.
    As $\mathcal{G}_T$ is finite this maximum always exists, however ties may occur making this a set-valued function. 
    Therefore a tie breaking procedure, determined on an application-specific basis, must be used to turn this into a function.
    The hyperparameter $\lambda >0$ determines the strength of the complexity penalty and can also be chosen on an application specific basis; $\lambda=1$ is a reasonable default.
    The scaling factor in terms of $|f|$ on the complexity term was included so that as the size $|f|$ of the sample set increases, the complexity penalty shrinks relative to the coverage score.



\section{Sample and Topic Metadata}
\label{s:metadata-migrations}

In addition to the hierarchical clustering itself, topic modelling entails generating metadata for each topic in the hierarchy.
Furthermore the samples themselves often come with metadata.
This metadata can be encoded in a relational database, but we still need to understand how to migrate data between these relational databases and the hierarchical topic model. 
 
\subsection{Sample Metadata Functor}
\label{ss:sample-metadata}
\begin{figure}
    \centering
    \begin{tabular}{c|c|c}
    Index & Author & PubYear \\
    \hline
     0 & Alice & 2014\\
     1 & Bob & 1998 \\
     ... & ... & ... \\
    \end{tabular}
    \caption{An example metadata table for our document set $X$, consisting of three columns.}
    \label{fig:metadata-table-example}
\end{figure}
Suppose that we have a table of metadata for our samples $X$. The relational model describes this information as a set of attribute sets $\{A_1,...,A_k\}$ and a relation $R_X\subseteq \prod_{i=1}^kA_i$ \cite{codd_relational_1970}. Each sample $x\in X$ defines a row of the table, which is an element $r_x = (r_1,...,r_k)\in R_X$. Let $A:=\prod_{i=1}^kA_i$ and let $\pi_i:A\to A_i$ denote the projection to the $i$th factor. We will assume that the first column is a primary key for $X$, meaning that the map $X\to A_1$ given by $x\mapsto \pi_1(r_x)$ is injective.

In categorical database theory, the attribute information is encoded as a schema category $M_X$. For a relational database with a single table having $k$ columns, such as our metadata table, the schema category is $M_X=\mathcal{P}([k])^{op}$. Here $\mathcal{P}$ denotes the power set, viewed as a poset category with arrows given by set inclusions.

Over this schema category, the attribute sets $A_i$ define a functor $\mathcal{A}:M_X\to \Set$ by
\begin{align*}
    \mathcal{A}(I) &= \prod_{i\in I} A_i,\\
    \mathcal{A}(I\to J) &= \pi_{I,J}, \qquad \forall: J\subseteq I,
\end{align*}
where $\pi_{I,J}:\mathcal{A}(I)\to \mathcal{A}(J)$ drops all the entries of $(a_i)_{i\in I}$ for $i\not\in J$.
The rows $R_X$ of a table with this schema are encoded as a subfunctor of $\mathcal{A}$, $D_X:M_X\to\Set$. This functor is defined by $D_X(I)= \{(r_i)_{i\in I}~:~ r \in R_X\}$ and $D_X( I\to J) = \pi_{I,J}|_{D_X(I)}$. 

For example, Figure \ref{fig:metadata-table-example} shows a potential metadata table containing three attributes, a primary key 'Index', the document's author 'Author', and its publication year 'PubYear'.  The schema category of this table is:
\[\begin{tikzcd}
	& {[\text{Documents]}} & \\
	{[\text{Index, Author}]} & {[\text{Index, PubYear}]} & {[\text{Author, PubYear}]} \\
	{\text{[Index]}} & {\text{[Author]}} & {[\text{PubYear}]} \\
	& \bot
	\arrow[from=1-2, to=2-1]
	\arrow[from=1-2, to=2-2]
	\arrow[from=1-2, to=2-3]
	\arrow[from=2-1, to=3-1]
	\arrow[from=2-1, to=3-2]
	\arrow[from=2-2, to=3-1]
	\arrow[from=2-2, to=3-3]
	\arrow[from=2-3, to=3-2]
	\arrow[from=2-3, to=3-3]
	\arrow[from=3-1, to=4-2]
	\arrow[from=3-2, to=4-2]
	\arrow[from=3-3, to=4-2]
\end{tikzcd}\]
Suppose the row of this example table corresponding to a sample $x$ is $r_x = (i_x, a_x, y_x)$, so that $R_X = \{r_x~:~x\in X\}$. 
Then 
\begin{align*}
    D_X(\text{Documents}) &= R_X, & D_X(\text{Index, Author}) = \{(i_x, a_x) :x\in X\},\\
    D_X(\text{Index}) &= \{(i_x):x\in X\}, &D_X(\text{Documents}\to \text{Index})(r_x) = (i_x),
\end{align*}
and so on. 
Using this framework, standard database operations become operations in the category of functors $M_X\to \Set$ \cite[\S 3.5]{fong_invitation_2019}. 

The schema category $M_X$ is determined by the columns of the table and the entries are determined by $D_X$ evaluated on the top object of the table: $D_X([k])$.
For this reason, it is sufficient to draw the schema category by only showing the top object of a table and its columns. 
Using this rule, the schema category of the table in Figure \ref{fig:metadata-table-example} would be depicted as: 
\[\begin{tikzcd}
	& {[\text{Documents]}} & \\
	{\text{[Index]}} & {\text{[Author]}} & {\text{[PubYear]}}
	\arrow[from=1-2, to=2-1]
	\arrow[from=1-2, to=2-2]
	\arrow[from=1-2, to=2-3]
\end{tikzcd}\]

Later we will consider subdatabases of $D_X$, meaning subfunctors $F\hookrightarrow D_X$. 
One important class of subdatabases are those generated by subsets of records $R_Y \subseteq R_X$ corresponding to subsets $Y\subseteq X$ of the samples.
For any categorical database $D:\mathcal{C}\to \Set$ let $\Sub(D)$ denote the category of subfunctors of $D$.
We define a functor $\mathcal{S}:\mathcal{P}(X) \to \Sub(D_X)$ by
\begin{equation}
    \label{eqn:def-of-S}
    \mathcal{S}(Y) := D_Y, \quad D_Y(J) := \{\pi_{[k],J}(r_x) ~|~ x\in Y\}, \quad \forall~ J \subseteq [k].
\end{equation}
If $Z\subseteq Y\subseteq X$, then for every $J\subseteq [k]$, $D_Z(J)\subseteq D_Y(J)$, and so $\mathcal{S}$ is indeed a functor.
Conversely, using $M_X$'s initial element $[k]$, we can define a functor $\Gamma:\Sub(D_X)\to \mathcal{P}(X)$ by $\Gamma(F) = \{x\in X ~:~ r_x \in F([k])\}$. 
This pair is an adjunction $\mathcal{S} \vdash \Gamma$ \cite{mac_lane_categories_1994}.

\subsection{Selection by a Match Condition}
\label{ss:selection-from-metadata}

To provide an example of how database operations can be described in categorical terms, we now review how to select data from $D_X$ according to a match condition. 
This is a two-step process. 
First, we compute a subset $\sigma\subseteq R_X$ of records that match the condition.
Then, we apply the functor $\mathcal{S}$ (Eqn. \ref{eqn:def-of-S}) to turn this subset into a categorical database.

Let $1$ denote the category with a single object $1$ and the identity morphism. 
A functor $F:1\to\Set$ is determined by the set $F(1)$, and this evaluation map is an isomorphism of categories $[1,\Set]\cong \Set$.
The category $1$ is the schema category for a table with a single column, and the set $F(1)$ contains the values in this column. A match condition (`where columns $I$ = value') for a selection query is encoded as a pair $(\bar{I}, v)$, where $\bar{I}:1\to M_X$ is the functor defined by $\bar{I}(1) = I$ and $v:\{\ast\}\to \mathcal{A}(I)$ is a map of sets sending $\ast$ to the desired value. 

The functor $\bar{I}$ between schema categories induces an adjoint pair of functors $\bar{I}^\ast \dashv \bar{I}_\ast$ between categories of databases:
\begin{itemize}
    \item $\bar{I}^\ast:[M_X,\Set]\to \Set$ mapping $F:M_X\to \Set$ to $F(I)$.
    \item $\bar{I}_\ast:\Set\to [M_X,\Set]$, taking the set $S$ to the functor defined by
    $$(\bar I_* S)(J)=\begin{cases}S,& I\subseteq J,\\ \{\mathrm{pt}\},& I\not\subseteq J.\end{cases}$$
\end{itemize} 
In the terminology of Spivak, $\bar{I}^\ast$ is the \emph{pullback data migration functor} and $\bar{I}_\ast$ is the \emph{right pushforward data migration functor}  \cite[\S 4]{spivak_functorial_2012}. Let $\eta:\mathrm{Id}_{[M_X,\Set]} \implies \bar{I}_\ast\circ \bar{I}^\ast$ denote the unit of this adjunction. 
There is a morphism $D_X\to \bar{I}_\ast \bar{I}^\ast \mathcal{A}$ given by composing the inclusion $D_X\to\mathcal{A}$ with $\eta(\mathcal{A})$.  
Let $\sigma(I,v)$ be the pullback of the diagram in $[M_X,\Set]$
$$D_X \xrightarrow{} \bar{I}_\ast\bar{I}^\ast \mathcal{A} \;\xleftarrow{\bar{I}_\ast v\,}\; \bar{I}_\ast\{\ast\}.$$
The set of records matching the condition $(I,v)$ is $\sigma(I,v)[k]$.

For example, let us describe the selection operation on the table in Figure \ref{fig:metadata-table-example} using the match condition `where Author="Alice"'. 
The match condition is expressed as a map $\text{Alice}:\ast \to \mathcal{A}(\text{Author})$ sending $\ast\to \text{"Alice"}$ and the result is defined by this pullback diagram:
\[\begin{tikzcd}
	{\sigma} && {\bar{I}_\ast\{\ast\}} \\
	{D_X} && {\bar{I}_\ast\bar{I}^\ast \mathcal{A}}
	\arrow[from=1-1, to=1-3]
	\arrow[from=1-1, to=2-1]
	\arrow["\ulcorner"{anchor=center, pos=0.125}, draw=none, from=1-1, to=2-3]
	\arrow["{\bar{I}_\ast \text{Alice}}", from=1-3, to=2-3]
	\arrow[from=2-1, to=2-3]
\end{tikzcd}\]
To obtain the set of documents matching this selection criterion we can expand the definition of the pullback:
\begin{align*}
    \sigma(\text{Documents}) &= \{(i,a,y)\in D_X(\text{Documents}): \pi_{\text{Author}}(i,a,y)=\text{Alice}(\ast)\} \\
    &= \{(i,a,y)\in R_X : a=\text{Alice}\}. 
\end{align*}
Finally, we apply $\mathcal{S}$ to obtain a subdatabase $D_\sigma = \mathcal{S}(\Gamma(\sigma))$ of $D_X$.
If we want the set of publication years of the documents authored by Alice, this would be $D_\sigma(\text{PubYear})$. 
The categorical description of other relational operations such as projections, joins and unions can be found in \cite{spivak_relational_2015}.

\subsection{Migrating Sample Metadata}
\label{ss:migrating-sample-metadata}

The metadata $M_X$ and cluster tree $T$ provide two different database schemata for organizing the samples $X$. 
To compose operations between these structures, such as "return the topic of the documents with author 'Bob'", we need a way to migrate data between these schemata. 

More precisely, we want to define a pair of functors between $\Sub(D_X)$ and $\Sub(\Clust)$. 
Here, $\Sub(\Clust)$ denotes the poset of monotone subfunctors of $\Clust$, meaning functors $G:\mathcal{G}_T\to [0,1]^X$ with $G(\phi)\sqsubseteq \Clust(\phi)$ for all $\phi \in \mathcal{G}_T$, with the pointwise order.
In Subsection \ref{ss:sample-metadata} we defined an adjoint pair $\mathcal{S}\vdash \Gamma$ between $\Sub(D_X)$ and $\mathcal{P}(X)$. 
Therefore it will now suffice to define a pair of functors between $\mathcal{P}(X)$ and $\Sub(\Clust)$.



Given a subset $Y\subseteq X$, let $\chi_Y:X\to \{0,1\}$ denote its characteristic function. This allows us to define a subfunctor $\Clust|_Y$ of the soft hierarchical clustering functor $\Clust:\mathcal{G}_T\to [0,1]^X$:
\begin{equation}
    \label{eqn:subfunctor-of-clust}
    (\Clust|_Y)(\psi) := (x \mapsto \chi_Y(x)\cdot\Clust(\psi)(x)), \qquad\forall \psi \in \mathcal{G}_T.
\end{equation}
Conversely, given a subfunctor $G\to \Clust$, for any $a\in [0,1]$ we can define a subset $Y_a(G)\subseteq X$ by  
\begin{equation*}
    Y_a(G) := \{x\in X ~:~ G(\top)(x)\geq a \}.
\end{equation*}
Putting this together, we have a family of functors for all $a\in [0,1]$:
\begin{align*}
    \Sigma_a:\Sub(\Clust) &\to \Sub(D_X), & 
    \Pi:\Sub(D_X)&\to \Sub(\Clust),\\
    G&\mapsto \mathcal{S}(Y_a(G)), & F&\mapsto \Clust|_{\Gamma(F)}.
\end{align*}

Using these functors, we can compose database operations defined using the schemata $M_X$ and $\mathcal{G}_T$.
For example, we can compute the query "return the topic of the documents with author 'Bob'" by first computing a selection subdatabase $F$ of $D_X$, then migrating it to $\Pi(F)$ and then computing $\Topic(\Pi(F)(\text{root}))$.

\begin{proposition}
    \label{prop:topic-preserves-specificity}
    Suppose $F_1\hookrightarrow F_2\hookrightarrow D_X$ are subfunctors of $D_X$. Then for every $\psi\in \mathcal{G}_T$
    \begin{equation*}
        \Topic(\Pi(F_1)(\psi)) \leq \Topic(\Pi(F_2)(\psi)).
    \end{equation*}
\end{proposition}
\begin{proof}
    Functors between posets are monotone maps, so $\Pi$ and $\Topic$ are both monotone. Thus $F_1\leq F_2$ implies $\Pi(F_1)(\psi)\sqsubseteq \Pi(F_2)(\psi)$, and monotonicity of $\Topic$ gives the final result.
\end{proof}
Proposition \ref{prop:topic-preserves-specificity} with $\psi=\top$ says that if you take a subset of samples in $X$ and apply a filter to obtain a smaller subset, the topic can only get more specific. 

\subsection{Migrating Cluster Metadata}
\label{ss:cluster-metadata}
In addition to sample metadata, we want to store metadata about the clusters, the vertices of $V(T)$.
As with the samples, this data can be organized into a table, modelled as a categorical database with a schema category $M_T$.
Suppose that there are $k'$ columns of metadata where the first column is a primary key for $V(T)$, meaning $M_T=\mathcal{P}([k'])^{op}$ and $V(T)$ injects into $A_1$.
Let $D_T:M_T\to \Set$ denote the database of metadata about the vertices of $T$. Let the set of records in the metadata table be $R_T = D_T([k'])$. The primary key induces an injection $\rho_T:V(T)\to R_T$. 

Similarly to the situation for document metadata, we want to define functors that migrate between databases over $\mathcal{G}_T$ and $M_T$.
Specifically, we'd like a functor $\Sub(D_T) \to \Sub(\Clust)$. 

Let $H$ be a subfunctor of $D_T$, corresponding to a subset of the clusters $V(T)$ arising as the results of queries on the metadata database $D_T$. Let $\chi_H:V(T)\to \{0,1\}$ be the characteristic function of the subset $\rho_T^{-1}\left(H([k'])\right) \subset V(T)$.
Then we define $F_H:\mathcal{G}_T\to [0,1]^X$ as
\begin{equation*}
    F_H(\psi) := \bigcup_{\substack{v\in V(T)\\v\leq \psi}} \chi_H(v)\cdot \Clust(v).
\end{equation*}
This is a pointwise colimit in $[0,1]^X$.
\begin{proposition}
    The map $F_H:\mathcal{G}_T\to [0,1]^X$ is a subfunctor of $\Clust$.
\end{proposition}
\begin{proof}
    We need to verify that $F_H$ is monotone and that $F_H(\theta)\sqsubseteq \Clust(\theta)$ for every $\theta \in \mathcal{G}_T$.
    For monotonicity, suppose $\psi\leq\phi$. Then
    $$
    {v:v\leq\psi,\ \chi_H(v)=1}
    \subseteq
    {v:v\leq\phi,\ \chi_H(v)=1}.
    $$
    Thus taking the union over this set we get $F_H(\psi) \leq F_H(\phi)$.
    Now suppose $v\leq \psi$ for some atomic formula $v\in V(T)$. Then $v\vee \psi = \psi$, and hence $\Clust(v) \sqsubseteq \Clust(\psi)$. 
    Then take the union over vertices $v\leq\psi$ on both sides to obtain:
    \begin{equation*}
       F_H(\psi) \sqsubseteq \Clust(\psi).
    \end{equation*}
\end{proof}
As both $\mathcal{G}_T$ and $[0,1]^X$ are lattices one may ask if $F_H$ is also a lattice homomorphism. In general, it only preserves joins and not necessarily meets.
For a counterexample consider two incomparable atomic formulas $a,b$ for which $\Clust(a)\cap\Clust(b)\neq 0$, such as Science and Math in the topic tree of Figure \ref{fig:small-topic-tree-example}. 
Suppose $H$ is a selection of vertices that includes both $a$ and $b$.
Then as there are no vertices below both $a$ and $b$, $F_H(a\wedge b)=0$, but on the other hand, $F_H(a)\wedge F_H(b) = \Clust(a)\cap \Clust(b) \neq 0$.

When computing $F_H(\psi)$ it suffices to take the union of $F_H(v)$ over atomic formulas $v\leq \psi$:
\begin{align}
    F_H(\psi) &= \bigcup_{\substack{v\in V(T)\\v\leq \psi}} \chi_H(v)\cdot \Clust(v)\notag\\
        &= \bigcup_{w\leq \psi} \bigcup_{v \leq w} \chi_H(v)\Clust(v)\notag\\
        &= \bigcup_{w\leq \psi} F_H(w). \label{eqn:F_H}
\end{align}
This is useful for computer implementation as it means we only need to store $F_H(v)$ for $v\in V(T)$ and can construct $F_H$ for general $\psi \in \mathcal{G}_T$ as needed. 

\section{Implementation of a Hierarchical Topic Database}
\label{s:implementation}

The mathematical description of a hierarchical topic database in Sections \ref{s:hierarchical-clustering-defn} and \ref{s:metadata-migrations} is abstract, but in fact the implementation of such a database and its operations reduces to simple matrix computations.
In this section we explicitly describe how one may store a hierarchical topic database in a computer, and how to compute Clust, RobustTopic and the metadata migrations.

\subsection{Data Model}
\label{ss:data-model}
    The initial data that we are indexing with a hierarchical topic database is a set $X=\{x_1,...,x_n\}$ of sample documents, such as text or images, with the following additional structure:
    \begin{enumerate}
        \item An embedding function $\text{Embed}:\{x_i\}\to \RR^d$, sending $x_i\mapsto p_i \in \RR^d$ for all $i=1,...,n.$
        \item A hierarchical topic model of the documents, meaning:
        \begin{enumerate}
            \item A rooted tree $T$ with vertices labelled $1,...,m$. Let $1$ denote the root of $T$.
            \item For each vertex $j \in V(T)$ a vector of cluster inclusion strengths\\ $\Clust(j):X\to [0,1]$, with $\Clust(1)(x)=1$ for all $x\in X$.
        \end{enumerate}
        \item A table of sample metadata $M_X$ with a primary key containing $X$ and with columns such as titles, authors, etc.
        \item A table of topic metadata $M_T$ with a primary key containing $V(T)$ and with columns such as keywords, summary phrases, etc.
    \end{enumerate}
    The hierarchical topic model can be stored using an $n\times m$ matrix $C$ valued in $[0,1]$, with $j$th-column given by the vector $\Clust(j)$ (so $ij$th entry $C_{ij} = \Clust(j)(x_i)$), and a set of tuples containing each vertex of the tree and its parent.
    For any topic $\phi \in \mathcal{G}_T$, the inclusion strength vector $\Clust(\phi)$ can be computed recursively from the rules
    \begin{itemize}
        \item $\Clust(j)$ is the $j$th column of $C$,
        \item $\Clust(\phi\wedge \psi) = \Clust(\phi)\cap \Clust(\psi)$, computed by taking the elementwise minimum,
        \item $\Clust(\phi\vee \psi) = \Clust(\phi)\cup \Clust(\psi)$, computed by taking the elementwise maximum.
    \end{itemize} 
    
    Hierarchical topic models obtained using soft clustering algorithms often result in every document having negligible but non-zero inclusion in every cluster. This is both meaningless in most applied contexts and inconvenient for memory usage in computation. 
    For this reason, we discretize the entries of the inclusion strength matrix $C$ to 8-bits, storing them as integers from $0-255$ and dividing by $255$ when we need values in $[0,1]$. 
    This will often make $C$ sparse, reducing memory usage and potentially improving the speed of matrix operations involving $C$.
    Of course, more or less precision could be used as appropriate for a given application.

    The documents and their metadata are stored in a standard table. The reference implementation uses Pandas. Similarly, the topic metadata for the vertices is stored in a table.
    The embedding vectors are stored as an array and we use PyNNDescent \cite{McInnes_PyNNDescent, dong_efficient_2011} for approximate nearest neighbour search.

    \subsection{Computing Metadata Selections}
    In Subsections \ref{ss:migrating-sample-metadata} and \ref{ss:cluster-metadata} we saw how metadata queries on the sample and topic metadata tables respectively can be used to obtain subfunctors of $\Clust$. 
    A subfunctor of $\Clust$ can have its value on vertices in the topic tree encoded as an $n\times m$ matrix $C'$ with $C'_{ij} \leq C_{ij}$ for every $i=1,...,n$ and $j=1,...,m$.
    If the subfunctor is also a lattice homomorphism, then this determines the subfunctor uniquely.
    Unfortunately, the subfunctor $F_H$ need not preserve meets. Nonetheless, we still only need to encode its value on the vertices of the cluster $T$ as equation \eqref{eqn:F_H} determines its value on a general element of $\mathcal{G}_T$ using only the vertex values.

    Selecting a subset of samples by metadata will reduce to multiplying $C$ by a diagonal matrix on the left, and selecting a subset of clusters will reduce to multiplying $C$ on the right, but using max-times matrix multiplication.

    A query on the sample metadata table returns a result subset $\sigma \subset X$ of samples. 
    Let $\chi_\sigma$ be the characteristic function of $\sigma$ and let $[\sigma]$ be the diagonal matrix defined by $[\sigma]_{ii} = \chi_\sigma(x_i)$. 
    Then the matrix $C_\sigma$ encoding the result subdatabase $\Clust|_\sigma$ is  
    $$C_\sigma = [\sigma]C.$$

    A query to the topic metadata table returns a subset $H \subset V(T)$ of topics. 
    Computing the matrix encoding the subdatabase $F_H$ is a bit more involved. 
    Let $[H]$ be the diagonal matrix with $[H]_{ii} = \chi_H(i)$ for $i=1,...,m$, let $A$ be the adjacency matrix of $T$, and let $\overline{A}$ be its transitive closure, meaning $\overline{A}_{ij} = 1$ if and only if $v_i\leq v_j$ in $T$.
    Then the matrix representing $F_H$ is 
    $$C_H = C \odot ([H]\bar{A}),$$
    where $\odot$ denotes the \emph{max-times} product, defined by $(M \odot N)_{ij} = \max_{k}(M_{ik}N_{kj})$.
    The $i$th column of this matrix stores the value of $F_H(i)$ for $i\in V(T)$. 
    To compute $F_H(\psi)$ for a general $\psi \in \mathcal{G}_T$, one must take the elementwise max over the columns $i \leq \psi$ of $C_H$. 


    \subsection{Incorporating Nearest Neighbour Selection}
    Allowing the user to pick out documents by embedding a search string and finding its $k$ nearest neighbour samples in the embedding space is an important utility for a practical implementation of topical search.

    Given a string $s$ let $B_k(s)$ denote the set of $k$-nearest-neighbours to $\text{Embed}(s)$ in $X$. Then the characteristic function of $B_k(s)$ is a function $\text{knn}(s):X\to\{0,1\}$. As with a metadata selection, we can apply the functor $-\cdot\Clust:[X,\{0,1\}]\to [T,[0,1]^X]$ to obtain a selection of the database $\Clust$. As a matrix computation, let $[\text{knn}(s)]$ be the diagonal $n\times n$ matrix with $[\text{knn}(s)]_{ii} = \text{knn}(s)(x_i)$, and 
    \begin{equation*}
        [\text{knn}(s)\cdot \Clust] = \text{[knn(s)]}C.
    \end{equation*}

    \subsection{Implementing RobustTopic}
    \label{ss:implementing-robust-topic}

    The robust topic function of Subsection \ref{ss:robust-topic} requires us to design algorithms to evaluate the surprise $S(f,\phi)$, the minimum description rank $\hat{r}(\phi)$ and to solve an optimization problem over the lattice $\mathcal{G}_T$.
    Evaluating the surprise is straightforward and computationally easy, but the other two tasks may be computationally difficult for large trees $T$.
    For this reason, I provide a simple greedy algorithm to compute RobustTopic. 
    If this model of hierarchical topic databases proves useful, then more research can be done on finding an efficient algorithm to compute RobustTopic more accurately. 
    
    Recall that for a fuzzy query set $f \in [0,1]^X$, we want to compute:
    \begin{equation*}
        \text{RobustTopic}(f) = \text{argmax}_{\phi\in\mathcal{G}_T} \left\{S(f,\phi) - \lambda \frac{\log(1+|f|)}{|f|} \hat{r}(\phi) \right\}, \qquad \lambda>0.
    \end{equation*}
    Fix $f$ and $\lambda$. For any geometric formula $\phi$, let 
    $$L(\phi) = S(f,\phi) - \lambda \frac{\log(1+|f|)}{|f|} r(\phi).$$
    Then our proposed greedy search algorithm is the following:
    \begin{algorithm}[Greedy RobustTopic algorithm]
        \label{alg:robusttopic-greedy}
    ~\\
    Input: A fuzzy set $f\in[0,1]^X$, a lattice $\mathcal{G}_T$ of topics, and a score function $L$.
    \begin{enumerate}[label=\arabic*., leftmargin=*]
        \item Initialize
        \[
            \phi_0 \in \arg\max_{v\in V(T)} L(v).
        \]

        \item For $i=0,1,\ldots$, repeat:
        \begin{enumerate}[label=(\alph*), leftmargin=*]
            \item \textbf{Disjunction step.}
            Choose
            \[
                v_i \in \arg\max_{v\in V(T)}
                    L(\phi_i\vee v).
            \]
            Set
            \[
                \psi_i =
                \begin{cases}
                    \phi_i\vee v_i,
                        & L(\phi_i\vee v_i)>L(\phi_i),\\
                    \phi_i, & \text{otherwise}.
                \end{cases}
            \]

            \item \textbf{Conjunction step.}
            Choose
            \[
                w_i \in \arg\max_{v\in V(T)}
                    L(\psi_i\wedge v).
            \]
            Set
            \[
                \phi_{i+1} =
                \begin{cases}
                    \psi_i\wedge w_i,
                        & L(\psi_i\wedge w_i)>L(\psi_i),\\
                    \psi_i, & \text{otherwise}.
                \end{cases}
            \]

            \item \textbf{Termination.}
            If neither step strictly improved the score,
            return $\phi_i$.
        \end{enumerate}
    \end{enumerate}
    \end{algorithm}
    In the reference implementation, Algorithm \ref{alg:robusttopic-greedy} is straightforwardly modified to beam search with default beam width of 3.
    In the optimization target $L(\phi)$ we are using the rank of the formula $r(\phi)$ in place of the true minimum description rank $\hat{r}([\phi])$ of the equivalence class of $\phi$. 
    Given a better algorithm $A$ for estimating the minimum description length of $[\phi]$ from $\phi$, one could improve Algorithm \ref{alg:robusttopic-greedy} by replacing $r(\phi)$ with $A(\phi)$.

    It is possible that the same equivalence class of formulas can appear multiple times in this search. For example, $[\psi] = [\psi \vee \phi \wedge \psi]$. 
    Fortunately, equivalent formulas have the same surprise score, so among the equivalent formulas considered during the search, the complexity penalty ensures the algorithm chooses the simpler presentation of the equivalence class. 
    In the next section we will present some results obtained by applying this algorithm to real-world data. 

\rowcolors{2}{white}{gray!15}
\section{Examples with Real-World Data}
\label{sec:examples}

    In this section we construct hierarchical topic models on two real-world text datasets and then demonstrate how one might use the database operations to do content recommendation and exploratory data analysis.
    In Subsection \ref{ss:20newsgroups}, we use a hard hierarchical topic model formed using human-generated topic labels, and in Subsection \ref{ss:ex-unsupervised}, we generate a soft hierarchical topic model using typical machine learning topic modelling tools.

    These examples are primarily meant to elucidate the semantics of the operations $\Clust$ and $\Topic$ as well as elements in the category $\mathcal{G}_T$ associated to a hierarchical topic model.
    In particular, we will see how conjunctions and disjunctions of primitive topics arise naturally in applications.
    These examples also demonstrate the utility and limitations of the RobustTopic function and Algorithm \ref{alg:robusttopic-greedy} for estimating it. 
    Every instance of a reported RobustTopic result in this section is an estimation computed using Algorithm \ref{alg:robusttopic-greedy}, modified to a beam search with width 3.

\subsection{Human-Curated Hierarchical Topic Model}
    \label{ss:20newsgroups}

    \begin{figure}
        \centering
        \includegraphics[width=0.75\linewidth]{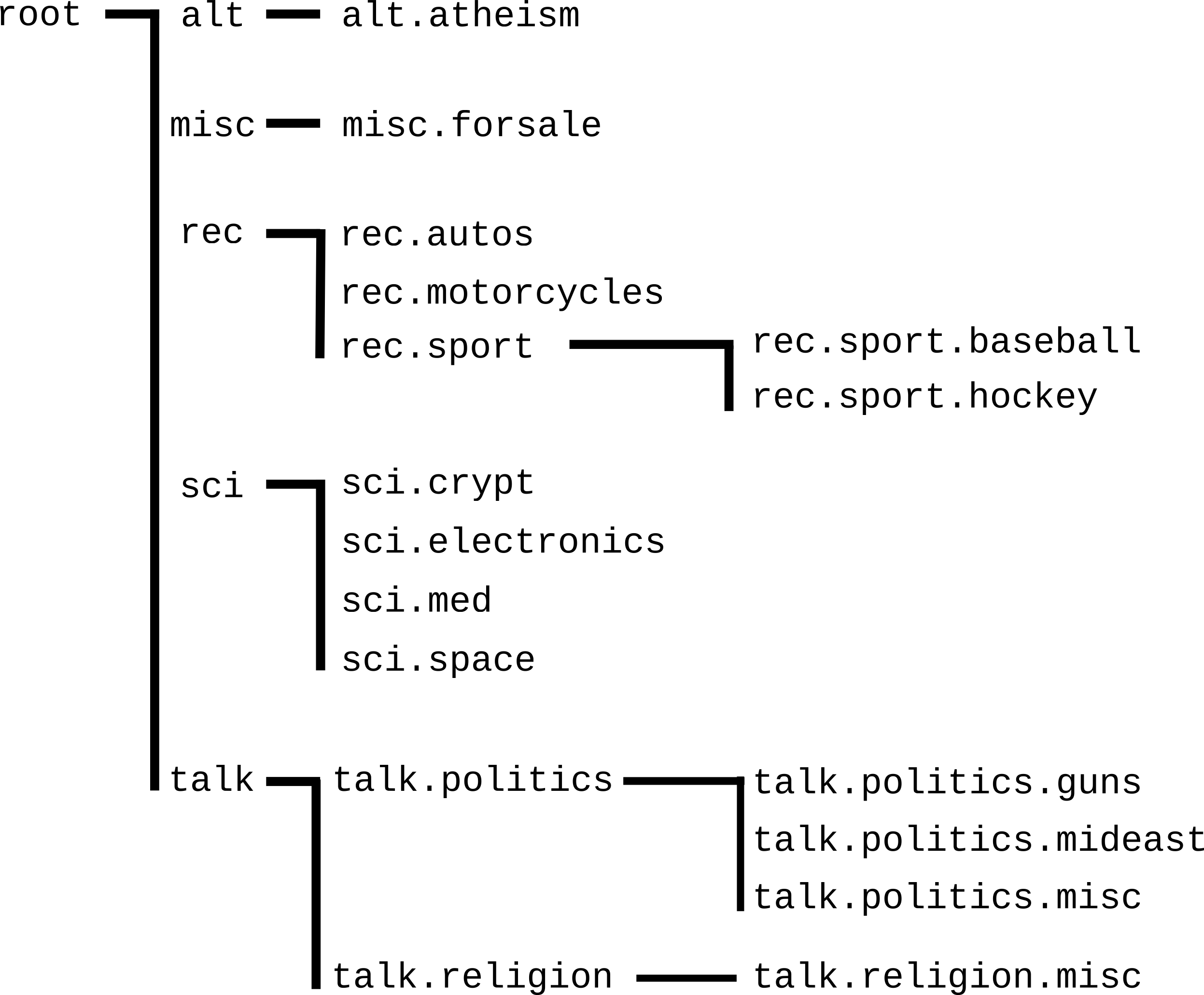}
        \caption{A subtree of the topic hierarchy formed by the newsgroups labels in the 20 Newsgroups dataset.}
        \label{fig:20ng-tree}
    \end{figure}

    \begin{figure}
        \centering
        \includegraphics[width=0.9\linewidth]{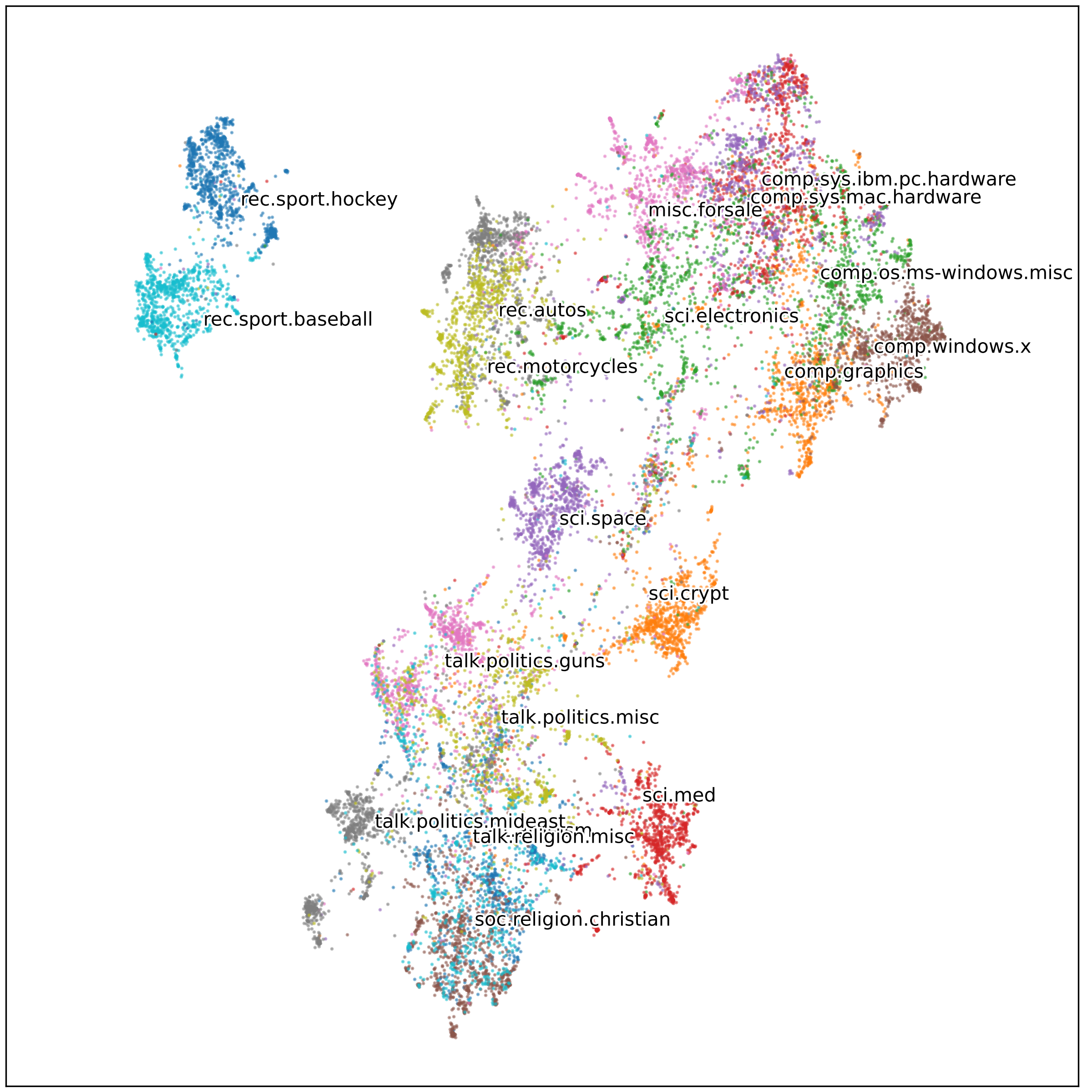}
        \caption{A scatter plot of the two-dimensional reduced embedding vectors of each post in the 20 Newsgroups dataset. Points are coloured by their newsgroup label. 
        The name of each newsgroup is plotted at the centroid of the corresponding post vectors.}
        \label{fig:20ng-umap}
    \end{figure}

    The 20 Newsgroups dataset is a common benchmark dataset for topic modelling tools that consists of  18,170 posts from Usenet discussion groups \cite{lang_newsweeder_1995}. Each post is labelled with a newsgroup, such as \verb|rec.sport.hockey| and these newsgroups are organized into a hierarchy (Figure \ref{fig:20ng-tree}). Thus 20 Newsgroups is equipped with a human-curated topic model that we can use for topical search. The newsgroups labels are binary, so this example does not use the full generality of the soft cluster tree; in the following subsection I will provide an example with $[0,1]$-valued inclusion strengths.
    Since each post belongs to exactly one leaf topic, the clusters of incomparable vertices in this hierarchy are disjoint. Consequently, the cluster strength vector of their conjunction is the zero vector.
    Therefore only disjunctions will appear (after simplifying conjunctions of comparable nodes) in geometric formulas obtained by applying RobustTopic to subsets of this database.

    For this experiment, I used a copy of the 20 Newsgroups dataset that has pre-computed embeddings into $\RR^{768}$ using the sentence-transformer model \verb|all-mpnet-base-v2| and reductions to two dimensions using UMAP \cite{mcinnes_umap_2018}. This pre-computed dataset is available on Hugging Face, at the URL:
    \begin{center}
        \href{https://huggingface.co/datasets/lmcinnes/20newsgroups\_embedded}{huggingface.co/datasets/lmcinnes/20newsgroups\_embedded}.
    \end{center}
    Instructions on how to reproduce this topic model using the Python reference implementation accompanying this paper are available at
    \begin{center}
        \href{https://thematic-search.readthedocs.io/en/latest/newsgroups.html}{thematic-search.readthedocs.io/en/latest/newsgroups.html}.
    \end{center}
    Figure \ref{fig:20ng-umap} is a scatter plot of the UMAP reduced vectors, coloured by newsgroup label. 

    In addition to the embeddings and cluster information, for each post $d_i$ we have two metadata items, \verb|post| which contains the text of the post and $\verb|newsgroup|$ which contains the name of the newsgroup to which it was posted. For each topic, we have only the metadata of its name, though one could compute additional information such as informative keywords.

    \subsubsection{RobustTopic outputs}

    Let $X$ be a set of samples constructed as follows. 
    Let $T$ be a subset of newsgroups tags in the topic tree, and for each $t \in T$ let us sample 100 posts from the dataset with topic $t$, without replacement.
    Then $X$ contains all these samples, plus $5*|T|$ additional noise samples chosen uniformly at random from the full dataset. For the remainder of this section we compute RobustTopic using Algorithm \ref{alg:robusttopic-greedy} with $\lambda=1$. 
    Here are examples of the computed robust topic for some choices of $T$:

    \begin{center}
    \begin{tabular}{|c|c|}
        Tags $(T)$ & Robust topic \\
        \hline
        \begin{tabular}{c} 
            \verb|rec.motorcycles|\\
        \end{tabular} & \verb|rec.motorcycles|
        \\
        \begin{tabular}{c} 
            \verb|rec.sport.hockey|\\
            \verb|rec.motorcycles|\\
        \end{tabular} & \verb|rec.sport.hockey| $\vee$ \verb|rec.motorcycles|
        \\
        \begin{tabular}{c} 
            \verb|rec.sport.baseball|\\
            \verb|rec.sport.hockey|\\
            \verb|rec.motorcycles|\\
        \end{tabular} & \verb|rec.sport| $\vee$ \verb|rec.motorcycles|
    \end{tabular}
    \end{center}

    In the first row we see that when we sample from a single topic (plus some noise), RobustTopic returns that topic. 
    In the second row we see that when we sample from two topics, the topic of the union is the disjunction of the two topics we sampled from.
    Finally in the third row, we sampled from the two subtopics of \verb|rec.sport|; instead of getting the disjunction of three topics, the two subtopics of \verb|rec.sport| are combined into a topic with lower complexity.  
    These examples also demonstrate the utility of passing from the topic tree to the Lindenbaum-Tarski algebra $\mathcal{G}_T$ as the second and third examples are best fit by topics that are disjunctions of the tree vertices.

    Here is another pair of outputs:

    \begin{center}
    \begin{tabular}{|c|c|}
        Tags $(T)$ & Robust topic\\
        \hline
        \begin{tabular}{c} 
            \verb|sci.med|\\
            \verb|sci.crypt|\\
        \end{tabular} & \verb|sci.crypt| $\vee$ \verb|sci.med| 
        \\
        \begin{tabular}{c} 
            \verb|sci.med|\\
            \verb|sci.crypt|\\
            \verb|sci.space|\\
        \end{tabular} & \verb|sci|
    \end{tabular}
    \end{center}

    When we sample from only two subtopics of \verb|sci|, we obtain a formula with a disjunction. 
    However when we sample from three of the four subtopics of \verb|sci|, RobustTopic identifies the topic as just \verb|sci|. The next table breaks down the score function (\ref{eqn:topic-score-function}) into the surprise and complexity terms for some candidate topics for this query:
    \begin{center}
    \begin{tabular}{|c|c|c|c|}
        Candidate topic & Score & Surprise & Complexity penalty\\
        \hline
        \verb|sci| & 1.47 & 1.47 & 0\\
        \verb|sci.crypt| $\vee$ \verb|sci.med| $\vee$ \verb|sci.space| & 1.70 & 1.73 & 0.037\\
        \verb|sci.crypt| $\vee$ \verb|sci.med| & 1.13 & 1.15 & 0.018
    \end{tabular}
    \end{center}
    From this table we see that the three-term disjunction is a better scoring topic, as it has higher surprise that outweighs the increase in complexity. 
    This demonstrates a limitation of the greedy algorithm; to find the three-term disjunction, the algorithm would have to first select one of three terms, then disjoin with a second term, then finally disjoin with the third term.
    However as we see in the third row of the table, these intermediate terms have a low score; thus Algorithm \ref{alg:robusttopic-greedy} (using a beam width of 3) converges to the local optimum \verb|sci|.

    \subsubsection{Querying for documents about science}
    The \verb|sci| branch of the 20 newsgroups topic tree contains topics about science. 
    Let's compare two ways that one may obtain a subset of samples that are about science:
    \begin{enumerate}
        \item Using $k$-nearest-neighbour semantic search on the string 'science' (or similar).
        \item Using the topic model to find posts under the \verb|sci| branch of the tree.
    \end{enumerate}
    Comparing the results of these methods is illustrative of the difference between the metric space topology induced on the dataset from the embedding model and the non-metric topology induced on the dataset from the hierarchical topic model. 
    The second operation is evaluating the functor $\Clust$ defining the hierarchical topic model at the topic $\phi = \verb|sci|$. 
    The formulas $\phi$ in the Lindenbaum-Tarski algebra are the summary objects that a human analyst naturally thinks about when working with topic models. The functor $\Clust$ interprets these formulas as fuzzy sets of samples.

\begin{figure}
        \centering
        \includegraphics[width=0.95\linewidth]{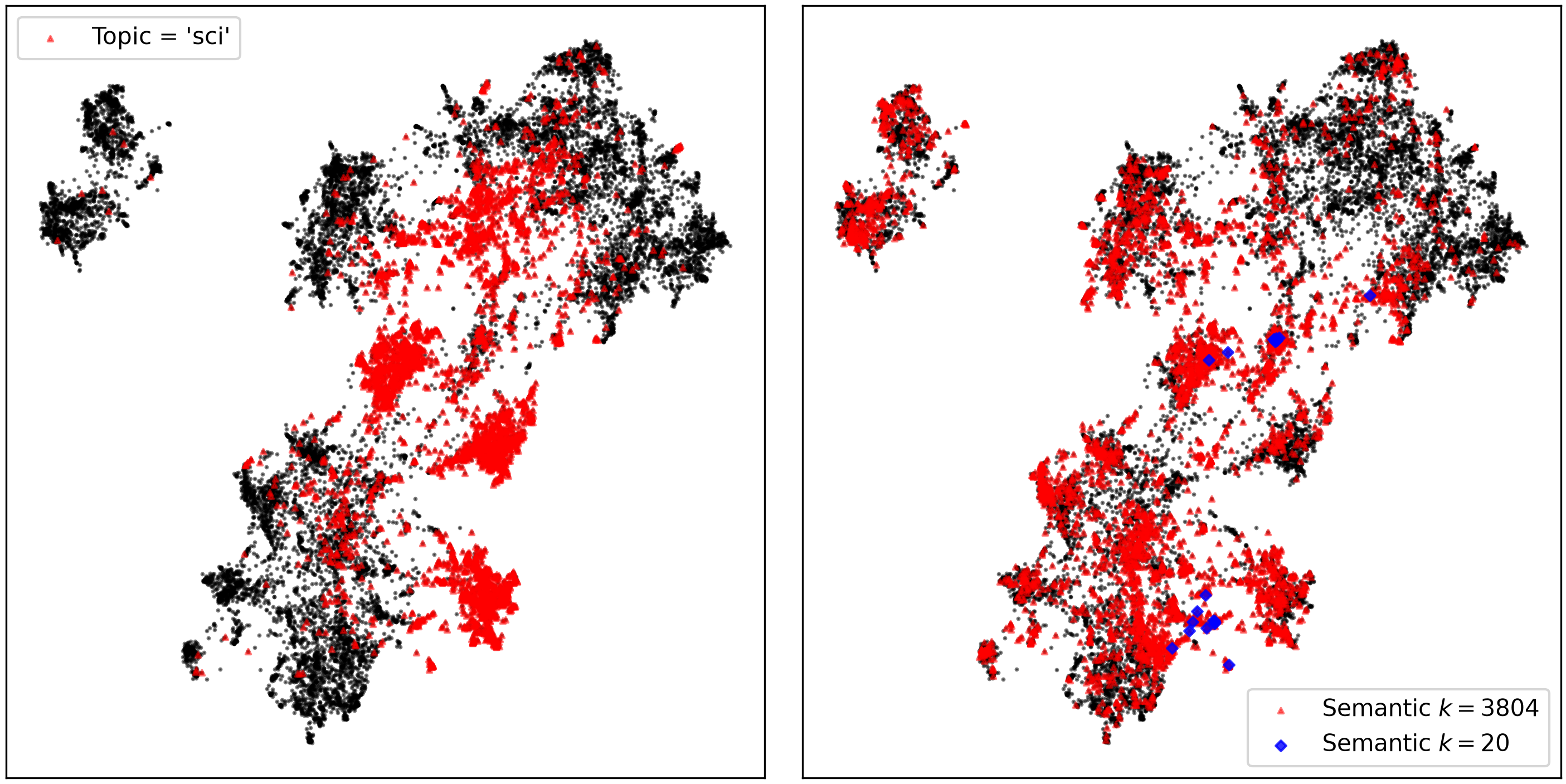}
        \caption{Comparison of the result set using selection by topic (left) and selection by semantic search (right) with $k=20$ and $k=3804$. Each point in the scatter plots is the dimensionally-reduced embedding of a single sample.}
        \label{fig:semantic-thematic-20ng}
    \end{figure}

    There are 3804 samples in the dataset labelled with the topic \verb|sci|, and thus 3804 results in the subset $\Clust(\verb|sci|)$. 
    Using semantic search with $k=20$, we obtain 20 posts that all match the topic of science well, but if we set $k=3804$ to get the same number of results as the topical filtering, then the relevance of many of the later posts to the topic of science is not high. For example:  
    \begin{pythonexample}{For the second straight game, California scored a ton of late runs to crush
the Brewhas. It was six runs in the 8th for a 12-5 win Monday and five in
the 8th and six in the 9th for a 12-2 win yesterday. Jamie Navarro pitched
seven strong innings, but Orosco, Austin, Manzanillo an ...}
    print(results.post.values[-1][0:280], "...")
    \end{pythonexample}
    Figure \ref{fig:semantic-thematic-20ng} provides a visual comparison of the result set using the UMAP reduced embedding vectors.  
    The result set obtained by topical filtering is considerably less dispersed in the UMAP representation than the result set obtained by nearest-neighbour search. Since the UMAP representation was constructed from semantic embeddings independently of the human-curated topic labels, this provides evidence that the topic hierarchy identifies subsets of documents that are coherent with respect to the semantic embedding geometry.
    Nearest-neighbour search selects documents according to their proximity to a single query embedding, while topical filtering selects documents according to membership in a hierarchically defined topic.
    For broad topics such as science, searching by topic can recover semantically coherent collections of documents that are not adequately represented by the nearest neighbours of a single query vector.

    \subsubsection{Semantic and topical integration}
    \label{sss:semantic-topical-integration}
    As an example of how the migration functors allow us to integrate semantic queries with the hierarchical topic information, consider the following content recommendation task.
    A new user of Usenet is interested in finding newsgroups that contain content they are interested in. One recommendation process could be
    \begin{enumerate}
        \item The user is prompted to enter a search string; a neural network embeds this into the same embedding space as the newsgroups post content.
        \item The set of $k$-nearest-neighbours to this embedding vector is computed.
        \item RobustTopic is used to compute a topic formula, a disjunction of newsgroups.
        \item The set of newsgroups appearing under the topics in the disjunction are\\ recommended to the user.
    \end{enumerate}
    If the recommendation task was to recommend \emph{posts} instead of newsgroups, we could then apply $\Clust$ to the topic from step 3 to obtain a set of posts to recommend. The resulting composition $\Clust\circ\operatorname{RobustTopic}$ is a robust analogue of the topical closure operator $\Clust\circ\Topic$ from Corollary \(\ref{cor:closure}\).

    The number $k$ of nearest neighbours can be used to influence the specificity of the results. As $k$ is increased, the resulting set of neighbours will be less semantically clustered and is therefore more likely to contain samples from multiple topics.
    For example, using our reference implementation, here are some search string inputs and outputs from Algorithm \ref{alg:robusttopic-greedy}:

    \begin{pythonexample}{comp.os.ms-windows.misc v comp.sys.ibm.pc.hardware v comp.graphics}
    search_string = "graphics drivers in windows"
    topic_db.q.neighbours(search_string, k=50).topic().formula
    \end{pythonexample}

    \begin{pythonexample}{soc.religion.christian v talk.religion.misc}
    search_string = "bible study"        
    topic_db.q.neighbours(search_string, k=100).topic().formula
    \end{pythonexample}

    \begin{pythonexample}{soc.religion.christian}
    search_string = "bible study"        
    topic_db.q.neighbours(search_string, k=50).topic().formula
    \end{pythonexample}
    
    \begin{pythonexample}{talk.politics.misc v rec.sport.baseball v talk.politics.guns}
    search_string = "ronald reagan versus jimmy carter"        
    topic_db.q.neighbours(search_string, k=100).topic().formula
    \end{pythonexample}
    The first three examples demonstrate that this process provides qualitatively plausible recommendations based on a search query. For the fourth example, it is not clear that \\
    \verb|rec.sport.baseball| is a good recommendation for the query.
    One limitation of this process is that when the query does not match any topic in the hierarchical topic model very closely, the resulting topic is essentially meaningless. For example:
    \begin{pythonexample}{rec.sport.hockey v rec.motorcycles v talk.politics.mideast}
    search_string = "ceci n'est pas un pipe"
    topic_db.q.neighbours(search_string, k=100).topic().formula
    \end{pythonexample}
    The topic score function (Equation \ref{eqn:topic-score-function}) lets us measure how well the nearest neighbours of a query align to a topic formula in the hierarchical topic model. 
    For example, the first three queries above have scores of 1.60, 2.06 and 2.17 respectively, whereas the final query has a score of 0.5. For completeness, the Reagan vs. Carter query has a score of 0.91. 
    Qualitatively, plausibility of the recommendation seems to increase as the score does, providing evidence that this score function is reasonable.

\subsection{Machine-Learned Hierarchical Topic Model}
    \label{ss:ex-unsupervised}

    In this section we will construct a hierarchical database using unsupervised machine learning tools.
    The dataset is English Wikipedia's \emph{vital level 5} articles, meaning the approximately 50,000 articles that Wikipedia's volunteer editors consider most important. 
    A snapshot of English Wikipedia articles was taken on 2026-07-01, then chunked and embedded using the Qwen3-Embedding-4B model by Sherlock Communications.
    This pre-embedded dataset is available on Hugging Face at the URL:
    \begin{center}
        \href{https://huggingface.co/datasets/Sherlock-Comms/wikipedia-en-2026-07-01-qwen3-embed-4b}{huggingface.co/datasets/Sherlock-Comms/wikipedia-en-2026-07-01-qwen3-embed-4b}.
    \end{center}
    After extracting the vectors for chunks coming from the vital level 5 articles, our dataset consists of 371,223 chunks coming from 43,115 articles. There were 6,641 vital level 5 articles that did not have corresponding vectors in the Sherlock Communications dataset, and are therefore missing from our model.
    Each chunk comes with the metadata of the title of the article and the section header from which it was taken, and the text of that chunk. 
    Let $Y=\{y_1,...,y_N\}$ denote the set of chunks and let $X=\{x_1,...,x_n\}$ denote the set of articles.

    \subsubsection{Topic model construction}
    To construct a hierarchical topic model for this dataset I used UMAP \cite{mcinnes_umap_2018} and the Toponymy Python package \cite{healy_toponymy_2024}.
    First, UMAP was used to reduce the 1024-dimensional embedding vectors from the dataset to 2 dimensions. 
    Then, the Toponymy library uses the embedding and reduced vectors along with the chunk text to construct a topic model as follows.
    First, a hierarchical clustering algorithm based on HDBSCAN is applied to the reduced vectors to obtain a tree $T$ and a binary cluster strength function $f_v:Y\to \{0,1\}$ for each cluster $v\in V(T)$.
    Then for each cluster Toponymy constructs metadata. First, informative keywords and example documents for each cluster are found with the BM25 algorithm. Then this metadata is passed to an LLM (Cohere Command-A-03-2025) that is asked to come up with a short descriptive title for the cluster.
    Finally the titles for each cluster undergo a disambiguation step, where the LLM is given titles that are similar and close in the cluster hierarchy and asked to modify the names to be more distinctive.
    An interactive visualization of the topic tree $T$ and an interactive scatter plot of the chunks coloured by topic are available online at the URL 
    \begin{center}
        \href{https://www.math.uwaterloo.ca/~kdruscit/wiki\_vital5\_datamap.html}{www.math.uwaterloo.ca/\~{}kdruscit/wiki\_vital5\_datamap.html}.
    \end{center}

    This process results in a hard hierarchical topic model for the set $Y$ of chunks. 
    The topic tree has 1445 vertices and depth 5. 
    Next I used this topic model to construct a soft hierarchical topic model for the set $X$ of articles.
    Many of the chunks $y\in Y$ are classified as noise, meaning that they are not included in any cluster. For these chunks, I assigned a soft inclusion strength vector by propagating the cluster labels from $y$'s $k$-nearest-neighbours in the embedding space. 
    This process might break the functoriality of the hierarchical clustering, so we then need to ensure functoriality by propagating inclusion strength up the tree. 
    Specifically, for every vertex $v$ we set $f_v(y) = \max_{w\leq v} f_w(y)$, computed in order from leaf nodes up to the root. 
    Note that $f_\top(y)=1$ for every chunk before label propagation and the label propagation does not change this property.

    To agglomerate the inclusion strength from chunks into an inclusion strength for articles, there are at least two reasonable approaches one might take.
    The first is to average the inclusion strength over an article's constituent chunks:
    $$\Clust^{\text{avg}}(v)(x) = \frac{1}{|\{y\in x\}|}\sum_{y\in x} f_v(y),$$
    where $y\in x$ means that $y$ is a chunk from article $x$.
    The second is to take the topic-wise maximum over an article's constituent chunks:
    $$\Clust^{\text{max}}(v)(x) = \max_{y\in x} f_v(y).$$
    The average, $\Clust^{\text{avg}}(v)(x)$, could be interpreted as "the proportion of chunks in the article $x$ that are about topic $v$". 
    In contrast,  $\Clust^{\text{max}}(v)(x)$ could be interpreted as "the extent to which any part of article $x$ is about topic $v$". 
    For this paper I use the $\Clust^{\text{max}}$ as the hierarchical topic model, because its semantics better match my intent on what the topic model should capture. 
    
    For example, consider the article \emph{Cambridge}, about the city in England. 
    This article has chunks with inclusion strength $1$ in the topic \emph{English Ceremonial Counties and Local Government}, and other chunks with inclusion strength $0$ in that same topic. 
    Under the maximum agglomeration scheme, \emph{Cambridge} has strength 1 in the topic \emph{English Ceremonial Counties and Local Government}, whereas under the average agglomeration scheme, the inclusion strength of this article in this topic depends on the length of the article and how much information is included about other aspects of the city of Cambridge. On the other hand, long articles have more opportunities for a chunk to be incorrectly classified and therefore under max agglomeration an article may be strongly included in more spurious topics, whereas under average agglomeration these incorrectly classified chunks only contribute weakly to the average.
    My intention is to capture a topical description of the city of Cambridge itself using the Wikipedia article as a proxy, and therefore the values for $\Clust$ should be robust under changes to the article text.
    If one instead was studying the Wikipedia article specifically, then it may be more sensible to use the average agglomeration scheme.

    As described in Subsection \ref{ss:data-model} the inclusion strengths are discretized to integers between 0 and 255. For all the computational examples in this section, we display the inclusion strengths scaled to $[0,255]$ instead of $[0,1]$ as this avoids printing long strings of decimals.

    \subsubsection{Querying by topic}

    Now we provide an example of querying this dataset using the hierarchical topic model.
    For example, suppose we are interested in North American geotopographical features. 
    Wikipedia maintains an article, \emph{Geology of North America}, indicating that this is a question of relevance.
    Instead, we can search for relevant articles using our hierarchical topic model. 
    There are primitive topics labelled 'Paleogeography, Tectonics, and Sediment Dynamics' (hereafter: Paleogeography) and 'North American Landscapes' (hereafter: N.A. Landscapes), and so their conjunction
    \begin{equation*}
        \text{Paleogeography} \wedge \text{N.A. Landscapes}
    \end{equation*}
    is a valid topic in our model, and we can use the cluster functor $\Clust$ to obtain a list of relevant articles, and because $\Clust$ is a lattice homomorphism we can compute $\Clust$ on this conjunction by taking the sample-wise minimum of $\Clust(\text{Paleogeography})$ and $\Clust(\text{N.A. Landscapes})$.
    Table \ref{tab:north-america-geology} shows the top 20 results by inclusion strength in this conjunction.

    \begin{table}[htbp]
        \centering
        \small
        \begin{tabular}{|p{4.2cm}rrrr|}
            \hline
            Article Title
            & \shortstack{Paleogeography}
            & \shortstack{N.A. Landscapes}
            & \shortstack{conjunction}
            & \shortstack{\# of\\Chunks} \\
            \hline
            Canadian Shield               & 255 & 255 & 255 &  8 \\
            Bay of Fundy                  & 255 & 255 & 255 &  5 \\
            Pacific temperate rainforests & 255 & 255 & 255 &  7 \\
            Ungava Peninsula              & 255 & 255 & 255 &  4 \\
            Backpacking (travel)          & 255 & 255 & 255 &  4 \\
            Antarctic ice sheet           & 255 & 255 & 255 &  7 \\
            Harbor                        & 255 & 255 & 255 &  6 \\
            Limestone                     & 255 & 255 & 255 & 17 \\
            Mountain                      & 255 & 255 & 255 & 10 \\
            Hudson Bay                    & 255 & 255 & 255 & 12 \\
            Esker                         & 255 & 255 & 255 &  6 \\
            Cenote                        & 255 & 255 & 255 &  8 \\
            Banc d'Arguin National Park    & 255 & 255 & 255 &  9 \\
            Tidal range                   & 255 & 255 & 255 &  3 \\
            Lake                          & 255 & 255 & 255 & 22 \\
            Geography of North America    & 255 & 255 & 255 & 13 \\
            Manicouagan Reservoir         &  74 & 255 &  74 &  4 \\
            Beringia                      &  71 & 255 &  71 & 10 \\
            Tectonic uplift               & 255 &  53 &  53 & 10 \\
            Lake Bonneville               &  40 & 255 &  40 & 12 \\
            \hline
        \end{tabular}
                \caption{Vital level 5 English Wikipedia articles sorted by their topical strength $\Clust(\text{Paleogeography} \wedge \text{N.A. Landscapes})$. Top 20 results are displayed.}
        \label{tab:north-america-geology}
    \end{table}

    In this table we see articles on named geographical features in North America, such as \emph{Canadian Shield} and \emph{Bay of Fundy}, as well as general geographical concepts present in North America such as \emph{Harbor} and \emph{Limestone}. However, \emph{Antarctic ice sheet} and \emph{Banc d'Arguin National Park} are not located in North America and \emph{Backpacking (travel)} is not a paleogeographic feature. Overall, this set of results is a qualitatively reasonable match for the query.

    \subsubsection{Topic of a sample subset}

    Given a set of samples from our dataset we can use the greedy RobustTopic algorithm \ref{alg:robusttopic-greedy} to find a topical description of the set.
    For example consider the sets
    \begin{align*}
        A &:= \{\text{Xhosa language}, \text{ Luganda}, \text{ Zulu language}, \text{ Shona language}\}\\
        B &:= \{\text{Southern Min}, \text{ Cantonese}, \text{ Hmong language}, \text{ Sichuanese dialects}\}.
    \end{align*}
    Set $A$ consists of languages native to central and southern Africa, and set $B$ consists of languages spoken by ethnic minorities in China. Applying Algorithm \ref{alg:robusttopic-greedy} we compute 
    \begin{align*}
        \text{RobustTopic}(A) &= ~\shortstack{\text{Grammatical Inflection}\\ \text{and Historical Linguistics}}~ \wedge \text{Southern and Central African Cities}.\\
        \text{RobustTopic}(B) &= \text{Chinese Language Varieties and Dialects}. 
    \end{align*}
    These results demonstrate one of the strengths of our categorical framework, in particular identifying topics of the model with the lattice of geometric formulas $\mathcal{G}_T$.
    In this dataset, there are enough articles about Chinese linguistics that the clustering algorithm identified them as a primitive topic. 
    However, the clustering algorithm did not identify a topic for African linguistics. 
    If we had defined the topics in the model to be the vertices of the topic tree $T$, then the best vertex to describe set $A$ would be 'Grammatical Inflection and Historical Linguistics'. 
    By permitting topics represented by geometric formulas, we are able to get a qualitatively more precise description of the set $A$.

    However, the effectiveness of the hierarchical topic model $\Clust$ is still limited by the quality of the topic tree $T$ and the hierarchical clustering. 
    For one example where our topic model underperforms expectations, consider the set 
    $$C := \{\text{Cambridge}, \text{ Oxford}\}.$$
    Cambridge and Oxford are both cities in England that are famous particularly for their world-renowned universities. Their Wikipedia articles both feature sections on their corresponding universities. For this reason, I expected that computing $\text{RobustTopic}(C)$ would return a result similar to (English Cities)$\wedge$(Universities). 
    Instead, we obtain
    \begin{align*}
        \text{RobustTopic}(C) =& ~\text{English Ceremonial Counties and Local Government}
    \end{align*}
    Indeed, there is no primitive topic in our model that directly "universities", so we could never see (English Cities)$\wedge$(Universities) as a topic, unless "universities" could be built up as a composite topic.
    As our topic model does not "see" the topic of universities, it is not able to find a concise topical description of set $C$ that distinguishes it from the larger set of English local governments.
    In this way, the descriptiveness of the set of geometric formulas $\mathcal{G}_T$ is still limited by the descriptiveness of the hierarchical topic tree $T$.

    \subsubsection{Exploratory data analysis}
    
    \begin{figure}
        \centering
        \includegraphics[width=0.3\textwidth]{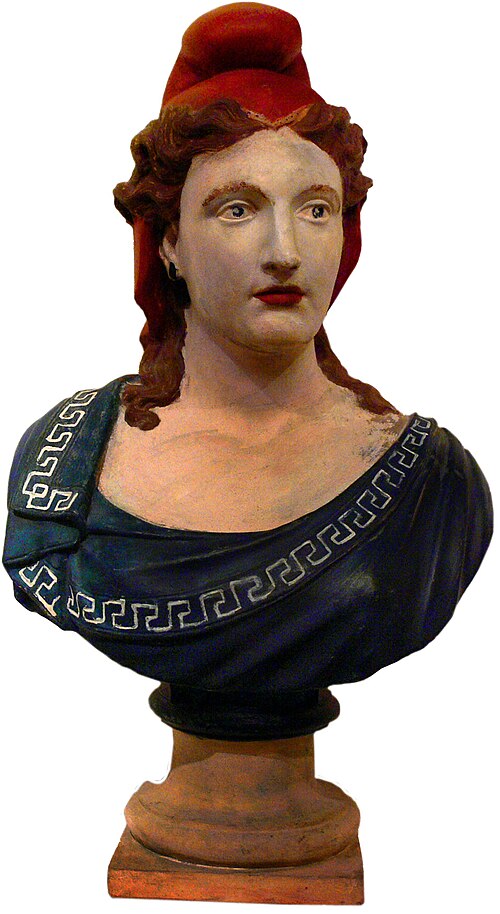}
        \caption{A bust of a woman wearing a Phrygian cap, located in the Palais du Luxembourg in Paris. Image retrieved from Wikipedia \cite{wikipedia-phrygian-cap}.}
    \end{figure}

    As one final application of this hierarchical topic model, we consider exploratory data analysis. 
    I was interested in finding articles that make connections between seemingly disparate topics, so I enumerated the pairs of incomparable vertices $v,v'$ in the topic tree and then looked for articles $x$ for which $\Clust(v \wedge v')(x)=255$. 
    The first example found is the article \emph{Phrygian cap}. 
    Quoting the Wikipedia article, the Phrygian cap is "a soft conical cap with the apex bent over, associated in antiquity with several peoples in Eastern Europe, Anatolia, and Asia" and "they came to signify freedom and the pursuit of liberty first in the American Revolution and then in the French Revolution, particularly as a symbol of Jacobinism (in which context it has been also called a Jacobin cap)" \cite{wikipedia-phrygian-cap}.
    This article was found in the conjunction of the topics 'Ancient Greek Tragedy and Comedy' and 'French Revolution: Key Figures and Events'. 

    Other interesting connections found this way were \emph{Orgasm} and \emph{Foreplay}, in the conjunction of 'Indian and Buddhist Philosophy' and 'Biological Sciences', and that there are 13 articles with strength 255 in the conjunction of 'Western Philosophy and Methodology' and 'Quantum Physics Foundations'. 
    These are all examples of  points of connection between fields of human knowledge that one can discover through the hierarchical topic modelling process described in this paper.

    \section{Future Research}

    There are multiple ways that this work can be improved or extended, both computationally and theoretically.

    Computationally, one could develop a better algorithm for computing the RobustTopic function of Section \ref{ss:robust-topic}, incorporating research on simplifying geometric formulas up to equivalence (e.g. \cite{couceiro:hal-01504010,willsey_egg_2021}) and efficient optimization over lattices.
    One could also study computational implementation of $\Clust$, $\Topic$ and the data migration functors, as well as methods to simplify composite queries for faster execution.
    A systematic empirical evaluation of the RobustTopic function would also be valuable. 
    This could include measuring its robustness under noise and changes in the topic model construction process, as well as human evaluation of the quality of the resulting topic formulas.  

    Theoretically, I believe the categorical database theoretic description of a hierarchical topic model could be extended more generally to other non-relational databases such as graph and vector databases, and this is a subject of ongoing research.
    Additionally, the category $\mathcal{G}_T$ is a \emph{frame}, as in pointless topology, and therefore there is a topological interpretation of the hierarchical topic model. 
    Developing this topological perspective further could help connect $\mathcal{G}_T$ to other topological data analysis methods, and may be a fruitful extension of this work.

    \section*{Acknowledgements}

    The author would like to thank Leland McInnes, John Healy and Aaron Smith for useful discussions during the course of this work. This work was supported by compute credits from a Cohere Labs research grant, and an NSERC PGS-D scholarship. 

\printbibliography
\end{document}

%% file: preamble.tex
\usepackage[utf8]{inputenc}
\usepackage{tgpagella,eulervm}
\usepackage[T1]{fontenc}
\usepackage{textcomp}
\usepackage{url}
 \usepackage{hyperref}
 \hypersetup{
     colorlinks,
     linkcolor={blue},
     citecolor={black},
     urlcolor={blue!80!black}
 }
\usepackage{graphicx}
\usepackage{float}
\usepackage{booktabs}
\usepackage{enumitem}
\usepackage{parskip}
\usepackage{emptypage}
\usepackage{subcaption}
\usepackage{multicol}
\usepackage[usenames,dvipsnames,table]{xcolor}
\usepackage[
	style=alphabetic,
	url=false,
	eprint=false,
    isbn=false,
]{biblatex}
\usepackage{xurl}
\usepackage{amsmath, amsfonts, mathtools, amsthm, amssymb}
\usepackage{mathrsfs}
\usepackage{cancel}
\usepackage{bm}
\usepackage[frozencache,cachedir=minted-cache]{minted2}
\setminted{fontsize=\small}

\newcommand\RR{\ensuremath{\mathbb{R}}}

\newcommand{\Sub}{\mathrm{Sub}}

\usepackage{systeme}
\let\svlim\lim\def\lim{\svlim\limits}
\let\implies\Rightarrow

\let\iff\Leftrightarrow
\let\epsilon\varepsilon
\usepackage{stmaryrd} 

\definecolor{correct}{HTML}{009900}

\newcommand\hide[1]{}

\usepackage{siunitx}
\usepackage{tikz}
\usepackage{tikz-cd}
\usepackage{quiver}
\usetikzlibrary{intersections, angles, quotes, calc, positioning}
\usetikzlibrary{arrows.meta}
\usepackage{pgfplots}
\pgfplotsset{compat=1.13}

\tikzset{
    force/.style={thick, {Circle[length=2pt]}-stealth, shorten <=-1pt}
}

\makeatother
\usepackage{thmtools}
\usepackage[framemethod=TikZ]{mdframed}
\mdfsetup{skipabove=1em,skipbelow=0em}

\usepackage[small]{titlesec}
\theoremstyle{definition}

\declaretheoremstyle[
    headfont=\sffamily\fontseries{bx}\selectfont\color{ForestGreen!70!black}, bodyfont=\normalfont,
    mdframed={
        linewidth=2pt,
        rightline=false, topline=false, bottomline=false,
        linecolor=ForestGreen, backgroundcolor=ForestGreen!5,
    }
]{thmgreenbox}

\declaretheoremstyle[
    headfont=\sffamily\fontseries{bx}\selectfont\color{NavyBlue!70!black}, bodyfont=\normalfont,
    mdframed={
        linewidth=2pt,
        rightline=false, topline=false, bottomline=false,
        linecolor=NavyBlue, backgroundcolor=NavyBlue!5,
    }
]{thmbluebox}

\declaretheoremstyle[
    headfont=\sffamily\fontseries{bx}\selectfont\color{NavyBlue!70!black}, bodyfont=\normalfont,
    mdframed={
        linewidth=2pt,
        rightline=false, topline=false, bottomline=false,
        linecolor=NavyBlue
    }
]{thmblueline}

\declaretheoremstyle[
    headfont=\sffamily\fontseries{bx}\selectfont\color{RawSienna!70!black}, bodyfont=\normalfont,
    mdframed={
        linewidth=2pt,
        rightline=false, topline=false, bottomline=false,
        linecolor=RawSienna, backgroundcolor=RawSienna!5,
    }
]{thmredbox}

\declaretheoremstyle[
    headfont=\sffamily\fontseries{bx}\selectfont\color{RawSienna!70!black}, bodyfont=\normalfont,
    numbered=no,
    mdframed={
        linewidth=2pt,
        rightline=false, topline=false, bottomline=false,
        linecolor=RawSienna, backgroundcolor=RawSienna!1,
    },
    qed=\qedsymbol
]{thmproofbox}

\declaretheoremstyle[
    headfont=\sffamily\fontseries{bx}\selectfont\color{NavyBlue!70!black}, bodyfont=\normalfont,
    numbered=no,
    mdframed={
        linewidth=2pt,
        rightline=false, topline=false, bottomline=false,
        linecolor=NavyBlue, backgroundcolor=NavyBlue!1,
    },
]{thmexplanationbox}

\numberwithin{equation}{section}

\declaretheorem[style=thmgreenbox, name=Definition]{definition}

\declaretheorem[style=thmbluebox, name=Algorithm]{algorithm}
\declaretheorem[style=thmredbox, name=Proposition]{proposition}
\declaretheorem[style=thmredbox, name=Theorem]{theorem}
\declaretheorem[style=thmredbox, name=Lemma]{lemma}
\declaretheorem[style=thmredbox, name=Corollary]{corollary}

\declaretheorem[style=thmproofbox, name=Proof]{replacementproof}
\renewenvironment{proof}[1][\proofname]{\vspace{-10pt}\begin{replacementproof}}{\end{replacementproof}}

\declaretheorem[style=thmexplanationbox, name=Proof]{tmpexplanation}

\usepackage{etoolbox}
\AtEndEnvironment{vb}{\null\hfill$\diamond$}%
\AtEndEnvironment{intermezzo}{\null\hfill$\diamond$}%

\makeatletter

\newcommand{\oefening}[1]{%
    \def\@oefening{#1}%
    \subsection*{Oefening #1}
}

\newcommand{\suboefening}[1]{%
    \subsubsection*{Oefening \@oefening.#1}
}

\newcommand{\exercise}[1]{%
    \def\@exercise{#1}%
    \subsection*{Exercise #1}
}

\newcommand{\subexercise}[1]{%
    \subsubsection*{Exercise \@exercise.#1}
}

\usepackage{xifthen}

\def\testdateparts#1{\dateparts#1\relax}
\def\dateparts#1 #2 #3 #4 #5\relax{
    \marginpar{\small\textsf{\mbox{#1 #2 #3 #5}}}
}

\def\@lesson{}%
\newcommand{\lesson}[3]{
    \ifthenelse{\isempty{#3}}{%
        \def\@lesson{Lecture #1}%
    }{%
        \def\@lesson{Lecture #1: #3}%
    }%
    \subsection*{\@lesson}
    \testdateparts{#2}
}

\makeatother

\usepackage{todonotes}
\usepackage[many,minted]{tcolorbox}
\DeclareTCBListing{pythonexample}{m}{%
    listing and comment,
    listing engine=minted,
    minted style=xcode,
    minted language=Python,
    minted options={linenos, numbersep=3mm},
    sharp corners,
    toprule=0pt, leftrule=0pt, rightrule=0pt,
    colframe=red!75!black,
    bicolor, colbacklower=white, colframe=black,
    fontlower=\ttfamily,
    comment={#1},
}
\tcbuselibrary{breakable}

\usepackage{import}
\usepackage{xifthen}
\usepackage{pdfpages}
\usepackage{transparent}

\author{Kaleb D Ruscitti}